\documentstyle[11pt]{article}
\def\div{{\,\rm div \,}}
\def\SOII{{\,\rm SO(2) \,}}
\def\SOIII{{\,\rm SO(3) \,}}
\def\adj{{\,\rm adj \,}}
\def\supp{{\,\rm supp \,}}
\def\B{{\,\bf B \,}}
\def\u{{\bf u}}
\def\v{{\bf v}}
\def\+M{{\,\rm M^{n\times n}_+ \,}}
\def\SO(n){{\,\rm SO(n) \,}}

\def\qfq{{\quad\mbox{for}\quad}}

\def\diag{{\,\rm diag \,}}
\def\supp{{\,\rm supp \,}}

\def\S{{\bf S}}

\def\lam{\lambda}

\def\X{{\cal X}}
\def\Q{{\cal Q}}

\def\A{{\cal A}}

\def\Q{{\cal Q}}

\newfont{\Blackboard}{msbm10 scaled 1200}

\newfont{\roma}{cmr10 scaled 1200}

\def\<{{\langle}}
\def\>{{\rangle}}

\def\Ga{\Gamma}
\def\var{\varphi}
\def\si{\sigma}

\def\a{\alpha}
\def\b{\beta}

\def\Om{\Omega}

\newtheorem{thm}{{}\hskip\parindent Theorem}[section]
\newtheorem{lem}{{}\hskip\parindent Lemma}[section]
\newtheorem{pro}{{}\hskip\parindent Proposition}[section]
\newtheorem{exl}{{}\hskip\parindent Example}[section]
\newtheorem{cor}{{}\hskip\parindent Corollary}[section]
\newtheorem{dfn}{{}\hskip\parindent Definition}[section]
\newtheorem{rem}{{}\hskip\parindent Remark}[section]

\def\dfrac{\displaystyle\frac}

\def\pl{\partial}
\def\rw{\rightarrow}

\def\na{\nabla}

\def\be{\begin{equation}}
\def\ee{\end{equation}}
\def\beq{\arraycolsep=1.5pt\begin{eqnarray}}
\def\eeq{\end{eqnarray}}

\def\R{I\!\!R}
\large
\title{Radial Deformations and Cavitation in Riemannian Manifolds with Applications to Membrane Shells}
\date{}
\author{
Peng-Fei YAO\\[0.3cm]
Key Laboratory of  Systems and Control\\
Institute of Systems Science,
Academy of Mathematics and Systems Science\\
Chinese Academy of Sciences, Beijing 100190, P. R. China\\
e-mail: pfyao@iss.ac.cn}

\begin{document}
\maketitle \footnote{This work is  supported by the National Science
Foundation of China, grants no. 60225003, no. 60334040, no.
60221301, no. 60774025, and no. 10831007. }
\begin{quote}
\begin{small}
{\bf Abstract} \,\,\,This study is a geometric version of Ball's work, Philos. Trans. Roy. Soc. London Ser. A 306 (1982), no.
1496, 557-611. Radial deformations in Riemannian manifolds are singular solutions to some nonlinear equations given by
constitutive functions and radial curvatures. A geodesic spherical cavity forms at the center of a geodesic ball in tension by means of given surface
tractions or displacements. The existence of such solutions depends on the growth properties of the constitutive functions and the radial curvatures.

Some close relationships are shown among radial curvature, the constitutive functions, and the behavior of bifurcation of a singular solution from a trivial solution.
 In the incompressible case the bifurcation depends on the local properties of  the radial curvature near the geodesic ball center but the bifurcation in compressible case
 is determined by the global properties of the radial curvatures.

A cavity forms at the center of a membrane shell of isotropic material placed in tension by means of given boundary tractions or displacements
when the Riemannian manifold under question is a surface of $\R^3$ with the induced metric. In addition, cavitation at the center of ellipsoids of $\R^n$ is also described if
the Riemannian manifold under question is $(\R^n g)$ where $g(x)$ are symmetric, positive matrices for $x\in\R^n.$
\\[3mm]
{\bf AMS(MOS) subject classifications}\,\,\,35D10
 \\[3mm]
\end{small}
\end{quote}

\tableofcontents

\setcounter{equation}{0}
\section{Introduction}
\def\theequation{1.\arabic{equation}}
\hskip\parindent The present paper is a geometric version of Ball's work \cite{Ball1}.

We investigate a class of singular solutions to the problems in which a hole forms in the center of a geodesic ball in a state of tension on a Riemannian manifold. This phenomenon of hole formation  is said to be {\it cavitation} by a terminology commonly used in the special case of an elastic fluid.

The study of cavitation  in the Euclidean space was initiated by Ball in the fundamental paper \cite{Ball1}.
The work of Ball was, in part, motivated by the work of \cite{GeLi} and subsequently developed by many authors
(see, e.g., \cite{Iw, Me, Si, Si1, SiSp, SiSp2, SiSp1,  St} and the review article \cite{HoPo}).

Let $(M,g)$ be a n-dimensional Riemannian manifold. Consider a body having strain energy $W$ in which the body occupies the open subset $\Om$ of $M.$ In a typical deformation in which a particle $x\in M$ is displaced to $\u(x)\in M$ the strain energy is given by
\be E(\u)=\int_\Om W(d\u)dg.\label{0.1}\ee The equilibrium equations of the body with zero body force are the Euler-Lagrange equations for the integral above. Solutions to these equilibrium equations are said to be equilibrium solutions.

Let $o\in M$ be given and let $\exp_o:$ $M_o\rw M$ be the exponential map. Radial deformations of a geodesic ball
$\Om=\{\,\exp_o\rho v\,|\,0\leq\rho<1,\,\,v\in M_o,\,\,|v|=1\,\}$ have the form
\be\u(x)=\exp_o\var(\rho)v\qfq x=\exp_o\rho v\in \Om.\label{0.2}\ee
Radial curvatures play a key role in the structure of radial deformations above. To get equilibrium equations for (\ref{0.2}) from (\ref{0.1}), we consider the case where radial curvatures on the geodesic sphere centered at $o$ with radius $t>0$ are the same, denoted by
$\kappa(t),$ for which we say $(M,g,o)$ is a {\it model} (\cite{GrWu}). We also assume that $W(F)$ can be expressed as a symmetric function $\Phi(v_1,\cdots,v_n)$ of the eigenvalues of $(F^TF)^{1/2}.$ For the incompressible case the only kinematically admissible deformations of the form (\ref{0.2}) are given by
\be \var(\rho)=\si^{-1}\Big(\si(\rho)+\si(A)\Big)\qfq \rho\geq0,\label{0.3}\ee where $A=\var(0)\geq0$ and  function $\si$ is defined by
\be\si(t)=\int_0^tf^{n-1}(s)ds\qfq t\geq0,\ee where $f$ is the solution to problem
\be f''(t)+\kappa(t)f(t)=0\qfq t>0;\quad f(0)=0,\quad f'(0)=1.\label{0.6}\ee For the compressible case $\var$ has to satisfy the radial equilibrium equation (Theorem \ref{nt1.1})
\be
[f^{n-1}(\rho)\Phi_1(\rho)]_\rho=(n-1)f^{n-2}(\rho)f'\circ\var(\rho)\Phi_2(\rho)\qfq x\in\Om,\quad
\rho(x)>0,\label{0.5}\ee where
$$\Phi_i(\rho)=\Phi_{v_i}\Big(\var'(\rho),\tau(\rho),\cdots,\tau(\rho)\Big),\quad\tau(\rho)
=\frac{f\circ\var(\rho)}{f(\rho)}\quad\qfq \rho>0,\quad i=1,\,\,2,$$ where $f$ is the solution to problem (\ref{0.6}).
Equation (\ref{0.5}) expresses the close relationship among radial curvature $\kappa,$ constitutive function $W$ and radial deformation $\var.$ Cavitation is equivalent to proving existence  of solutions to problem (\ref{0.5}) such that $\var(0)>0.$

Consider existence of equilibrium solutions in the incompressible case. The assumption on the radial curvature such that formula (\ref{0.3}) makes sense is the following
\be \int_0^{\max\{1\var(1)\}}s\kappa_+(s)ds\leq1,\label{0.7}\ee where $\kappa_+(s)=\max\{0,\kappa(s)\}.$ We show (Theorem \ref{t2}) that under assumption (\ref{0.7}) if $A>0$ then (\ref{0.3}) generates an equilibrium solution if and only if
\be\frac{v^{n-1}}{(v^n-1)^2}\hat{\Phi}'(v)\in L^1(\delta,\infty)\qfq\delta>1,\label{0.8}\ee where $\hat{\Phi}(v)=\Phi(v^{n-1},v,\cdots,v).$ Conditions (\ref{0.8}) are the same as in \cite{Ball1} in the Euclidean space. Let $P$ be the radial component of the Piola-Kirchhoff stress at $\rho=1$ and let $T$ be the Cauchy stress. Under assumptions (\ref{0.7}) and (\ref{0.8}) and with the choice $T(0)=0$ a critical value $P_{cr}$ of P for solutions of bifurcation is given by
\be P_{cr}=\int_1^\infty\frac{1}{v^n-1}\hat{\Phi}'(v)dv,\label{0.10}\ee which is again the same as in \cite{Ball1} in the Euclidean space.

In the compressible case the establishment of existence of cavitating equilibrium solutions is much more complicated. We use some similar assumptions on the growth properties of the constitutive function $W$ as in \cite{Ball1} to analyze  equilibrium solutions.  The displacement boundary value problem in which $\var(1)=\lam>0$ is specified is concerned.
A solution to problem (\ref{0.5}) is said to be {\it regular} if $\var(0)=0.$ Since there are no explicit formulas for
regular equilibrium solutions in general, we establish some estimates of regular equilibrium solutions from below and above (Theorems \ref{t3.2} and \ref{t3.3}) under the radial curvature assumption (\ref{0.7}). Using these estimates for regular equilibrium solutions and under the  radial curvature assumptions
\be\int_0^\infty s\kappa_+(s)ds\leq1,\quad\int_0^\infty s\kappa_-(s)ds<\infty, \label{0.11}\ee where $\kappa_-(s)=\max\{0,-\kappa\},$ we show (Theorem \ref{t4.4}) that for $\lam$ large enough there is a unique radial minimizer $\var$ of $E$ with $\var(1)=\lam$ and $\var(0)>0,$ which is a stable cavitating equilibrium solution.

One of the direct applications of the analysis here is the cavitation problem of membrane shells. Let $M$ be a surface in $\R^3$
 with the induced metric $g.$ Suppose that the middle surface of a shell is a bounded open set $\Om\subset M.$ First we show (Proposition \ref{p6.1}) that the $\Ga$-limit membrane shell, given in \cite{LeRa}, takes the form (\ref{0.1}) if all deformations of the middle are confined in $M.$ So we assume that the  membrane shells have their stored energies in the form (\ref{0.1}) to study their cavitation problems. In particular, the following surfaces of revolution are concerned:
 $$M=\{\,(x,\psi(r))\in\R^3\,|\,x=(x_1,x_2)\in\R^2,\,\,r=|x|\,\},$$ where $\psi$ is a $C^2$ function on $[0,\infty).$
 Then $(M,g,o)$ is a model where $g$ is the induced metric of $M$ from $\R^3$ and $o=(0,0,\psi(0)).$  The radial curvature is given by
 $$\kappa(t)=\frac{\psi'(\zeta(t))\psi''(\zeta(t))}{\zeta(t)(1+\psi'^2(\zeta(t)))^2}\qfq t\geq0,$$ where function $\zeta(t)$ is defined by equation
 $$ t=\int_0^{\zeta(t)}\sqrt{1+\psi'^2(s)}ds\qfq t\geq0.$$ In the incompressible case a critical value of the radial component of the Piola-Kirchhoff stress $P$ at $\rho=1$  for solutions of bifurcation is given by (\ref{0.10}) under the assumptions (\ref{0.7}), (\ref{0.8}) and $T(0)=0.$ In the compressible case when $\var(1)=\lam$ is large enough the stable solution is cavitating ($\var(0)>0$) under the assumptions (\ref{0.11}), $T(0)=0$ and the growth assumptions of $W.$

 Consider the Riemannian manifold $(\R^n,g)$ where $g=G(x)$ are  symmetric and positive matrices  for $x\in\R^n.$
 The radial deformation theory of Sections 2-4 describes that a  a hole forms in the center of an  ellipsoid in a state of tension in Section 6.

\setcounter{equation}{0}
\section{Equilibrium Equations for Radial Deformations on a Model}
\def\theequation{2.\arabic{equation}}
\hskip\parindent We make some preparations for our problems.
Let $(M,g)$ be a $n$-dimensional Riemannian
manifold with an orientation  and let $\Om\subset M$ be an open
set. A map $\u:$ $\Om\rw M$ is said to be a {\it deformation.} Let $\u:$
$\Om\rw M$ be a deformation. We define {\it the deformation gradient}
$d\u$ of $\u$ as a bilinear functional on $M_{\u(x)}\times M_x$ by
\be d\u(Y,X)=\<Y,\u_*X\>\circ\u(x)\qfq Y\in M_{\u(x)},\quad X\in
M_x,\quad x\in \Om,\label{1.3}\ee where $\<\cdot,\cdot\>=g$ is the
Riemannian metric.

Let $W$: $M_+^{n\times n}\rw \R$ be a function where $M_+^{n\times
n}$ is the set of real $n\times n$ matrices with positive
determinant.   We denote by $\SO(n)$ the special orthogonal group
on $\R^n$. We assume that \be W(F)=W(QF)=W(FQ)\qfq F\in
M_+^{n\times n},\quad Q\in\SO(n).\label{2}\ee

Let $x\in M$ be given. Let $\{e_i\}$ and $\{E_i\}$ be orthonormal
bases of $(M_x,g(x))$ and $(M_{\u(x)},g\circ\u(x))$ with the
positive orientation, respectively. We define \be
W(d\u)=W(F),\label{3}\ee where
\be F=\Big(d\u(E_i,e_j)\Big).\label{4*}\ee We have the following.

\begin{lem}\label{l1.1} Let $W$ satisfy  $(
\ref{2})$. Then the definition of $(\ref{3})$ is independent of
the selections of $\{e_i\}$ and $\{E_i\}$.
\end{lem}

{\bf Proof}\,\,\,Let $\{\hat{e}_i\}$ and $\{\hat{E}_i\}$ be the
different selections. Let $$\hat{e}_i=\sum_{j=1}^nq_{ij}e_j,\quad
\hat{E}_i=\sum_{j=1}^nr_{ij}E_j.$$ Then $$
\<\hat{E}_i,\u_*\hat{e}_j\>\circ\u(x)=\sum_{kl}r_{ik}\<E_k,\u_*e_l\>q_{jl},$$
that is,
$$\Big( d\u(\hat{E}_i, \hat{e}_j)\Big)=R\Big(
d\u(E_i,e_j)\Big)Q^T,$$ where $R=\Big(r_{ij}\Big)$ and
$Q=\Big(q_{ij}\Big)$ are in $\SO(n)$. Then the lemma follows from
(\ref{2}).    \hfill$\Box$

\begin{rem}Let $M=\R^n$ with the Euclidean metric and
$\u=(u_1,\cdots,u_n).$ Then $$\<\frac{\pl}{\pl
x_i},\u_*\frac{\pl}{\pl x_j}\>=\frac{\pl u_i}{\pl x_j}.$$
\end{rem}

In a typical deformation in which the point $x\in \Om$ is
displaced to $\u(x)\in M$ the  energy of $\u$ is given by \be
E(\u)=\int_\Om W(d\u)dg,\label{4}\ee where $dg$ is the volume
element of $M$ in the metric $g$.

\begin{pro}\label{pn2.1} Let $x=(x_1,\cdots,x_n)$ be a local coordinate system on $(M,g)$
and $$g=\sum_{ij=1}^ng_{ij}(x)dx_idx_j.$$ Let
$G=\Big(g_{ij}\Big).$ Then the equilibrium equations are the
Euler-Lagrange equations for $(\ref{4})$ \be
\sum_{lij}\frac{\pl}{\pl x_l}\Big(\frac{\pl W(d\u(x))}{\pl
F_{ij}}\a_{ip}(\u(x))\a^{jl}(x)\sqrt{\det G(x)}\Big)=0\qfq 1\leq
p\leq n,\label{5}\ee where
$$\Big(\a_{ij}\Big)=G^{1/2},\quad \Big(\a^{ij}\Big)=G^{-1/2}.$$
\end{pro}

{\bf Proof}\,\,\, Let $\v:$ $M\rw M$ be a deformation and let
$$e_i(x)=\sum_{j=1}^n\a^{ij}(x)\pl x_j\qfq1\leq i\leq n.$$ Then $e_1,$ $\cdots,$ $e_n$ form an orthonormal basis of $M_x$
and $E_1,$ $\cdots,$ $E_n$ is an orthnonrmal basis of $M_{\v(x)},$
where $E_i=e_i(\v(x))$ for $1\leq i\leq n.$ Moreover, it follows
that \beq\<E_i,\v_*e_j\>&&= \sum_{kl}\a^{ik}(\v(x))\<\pl
x_k,\v_*\pl
x_l\>\a^{jl}(x)\nonumber\\
&&=\sum_{klh}\a^{ik}(\v(x))g_{kh}(\v(x))\a^{jl}(x)\frac{\pl
v_h}{\pl
x_l}\nonumber\\
&&=\sum_{lh}\a_{ih}(\v(x))\a^{jl}(x)\frac{\pl v_h(x)}{\pl
x_l}.\nonumber\eeq

Let
$$I(\varepsilon)=E(\u+\varepsilon\v).$$
We have \beq I'(0)&&=\sum_{ij}\int_\Om\frac{\pl W}{\pl
F_{ij}}\<E_i,\v_*e_i\>\sqrt{\det G(x)}dx\nonumber\\
&&=\sum_{ijlh}\int_\Om\Big\{\frac{\pl}{\pl x_l}\Big(\frac{\pl
W}{\pl F_{ij}}\a_{ih}(\v(x))\a^{jl}(x)\sqrt{\det
G(x)}v_h\Big)\nonumber\\
&&\quad-\Big[\frac{\pl}{\pl x_l}\Big(\frac{\pl W}{\pl
F_{ij}}\a_{ih}(\v(x))\a^{jl}(x)\sqrt{\det
G(x)}\Big]v_h\Big\}dx.\nonumber\eeq Equations (\ref{5})
follow. \hfill$\Box$

\begin{rem}Let $M=\R^n$ with the Euclidean metric. Equations $(\ref{5})$ become
$$\sum_{j=1}^n\frac{\pl }{\pl x_j}\Big[\frac{\pl W(d\u)}{\pl F_{ij}}\Big]=0\qfq
1\leq i\leq n.$$
\end{rem}

\begin{rem} For a general deformation $\u,$ the problem
$(\ref{5})$ may be very complicated. Let
$$W(F)=\frac{1}{2}|F|^2\qfq F\in M_+^{n\times n}.$$ Then
equations $(\ref{5})$ are
$$\Delta u_p+\sum_{ijkl=1}^ng^{ij}(x)\Ga_{kl}^p(\u(x))\frac{\pl
u_k}{\pl x_i}\frac{\pl u_l}{\pl x_j}=0\qfq1\leq p\leq n,$$ where
$\Delta$ is the Laplacian on $M$ in the metric $g$ and
$\Ga_{ij}^k$ are the Christoffel symbols. Solutions of the above
equations are called harmonic maps, see Lemma $8.1.1$ in
$\cite{Jo}.$
\end{rem}

We are interested in radical deformations that are introduced below.

Let $o\in M$ be fixed and let $\exp_o:$ $M_o\rw M$ be  the
exponential map at the point $o$ in the metric $g.$ For any $x\in
M,$ there exists a pair $(\rho, v)$ with $\rho\geq0$  and
 such that\be
 x = exp_o\rho v \label{6}\ee where $v=v(x)\in S_o$ and
 $S_o$
is the unit sphere of $(M_o,g(o)).$  Let $d(x, y)$ be the distance
function from $x$ to $y$ in the metric $g.$ Then $\rho= d(o, x).$

\begin{dfn} A map $\u:$ $M\rw M$ is said to be a radical
deformation with respect to $o\in M$ if there is a function
$\var:$ $[0,\infty)\rw\R$ such that
\be\u(x)=\exp_o\var(\rho) v\qfq x=\exp_o\rho v\in
M.\label{7}\ee
\end{dfn}

We shall solve the problem (\ref{5}) when $\u$ is a radical
deformations under appropriate assumptions on the constitutive
function $W$ and on the geometric properties of the metric $g.$ For this end, we need to computer $d\u$ first.

Let $\X(M)$ be all vector fields on $M.$ Let D be the Levi-Civita
connection of the metric $g.$ Let X and Y be vector fields on $M.$
The curvature operator is a map ${\bf R}_{XY} : \X(M)\rw\X(M)$,
given by $${\bf R}_{XY} Z =-D_XD_Y Z + D_Y D_XZ + D_{[X,Y ]}Z$$
for all $Z\in\X(M),$ where $[\cdot,\cdot]$ is the Lie bracket
product. Let $\gamma(t)$ be a geodesic with $|\dot{\gamma}(t)|=1$ initiating from the point $o.$
A vector field $J:$ $[0,\infty)\rw M_{\gamma}$ is called a Jacobi
field along $\gamma$ if
$$\ddot{J}(t)+{\bf R}_{\dot{\gamma}(t)J}\dot{\gamma}=0\qfq
t\geq0.$$ In addition, a Jacobi field $J$ is said to be normal if
$$\<J(t),\dot{\gamma}(t)\>=0\qfq t\geq0.$$

\begin{pro}\label{pl3} Let $\u$ be a radical deformation and let
$\rho=\rho(x)$ be the distance function in the metric $g$ from
$x\in M$ to $o.$ Then \be
\u_*D\rho(x)=\var'(\rho)D\rho(\var(\rho))\qfq x\in M,\quad
x\not=o.\label{9n}\ee

Let $J(t)$ be a normal Jacobi field along the geodesic $\gamma(t)$
with $J(0)=0.$ Then \be \u_*J(t)=J(\var(t))\qfq t\in\R.\ee
\end{pro}

{\bf Proof}\,\,\,Formula (\ref{9n}) follows from expression
(\ref{7}).

Let $\gamma(t)=\exp_otv$ where $v\in M_o$ with $|v|=1.$ Since $\<v,\dot{J}(0)\>=0,$ there is a curve $\si:$ $[0,1]\rw M_o$
such that
$$|\si(s)|=1,\quad\si(0)=v,\quad \dot{\si}(0)=\dot{J}(0).$$
Let $$\a(t,s)=\exp_ot\si(s)\qfq (t,s)\in[0,\infty)\times[0,1].$$
Then $$J(t)=\a_s(t,0)=t\exp_o \dot{J}(0)\qfq t\in[0,\infty).$$

By definition, we have
$$\u(\a(t,s))=\exp_o\var(t)\si(s)\qfq (t,s)\in[0,\infty)\times[0,1],$$
which yields
$$\u_*J(t)=\var(t)\exp_{o*}\dot{J}(0)=J(\var(t)).$$  \hfill$\Box$\\

For any $v\in M_o$ with
$|v|=1,$  there is a unique $t_0(v)>0$ (or $t_0(v)=\infty$) such
that the normal geodesic $\gamma(t)=\exp_otv$ is the shortest on the
interval $[0,t_0).$  Let
$$C(o)=\{\,\,t_0(v)v\,\,|\,\,v\in M_o,\,\,|v|=1\,\,\},\quad
\Sigma(o)=\{\,\,tv\,\,|\,\,v\in M_o,\,\,|v|=1,\,\,0\leq t<
t_0(v)\,\,\}.$$ The set $\exp_oC(o)\subset M$ is said to be the
cut locus of $o$ and the set $\exp_o\Sigma(o)\subset M$ is called
the interior of the cut locus of $o.$ Then
$$M=\exp_o\Sigma(o)\cap\exp_oC(o).$$ Furthermore, $\exp_o:$ $\Sigma(o)\rightarrow
\exp_o\Sigma(o)$ is a diffeomorphism and $C(o)$ is a zero measure
set on $M_o$. Then $\exp_o C(o)$ is a zero measure set on $M$
since it is the image of the zero measure set $C(o)$, that is,
$\exp_o\Sigma(o)$ is $M$ minus a zero measure set.

Let $\psi:$ $M_o\rw M_o$ be a linear opeator.  We define a map $\Psi:$ $\exp_o\Sigma(o)\rw\exp_o\Sigma(o)$ by
\be\Psi(x)=\exp_o\rho\psi v\qfq x=\exp_o\rho v\in \exp_o\Sigma(o).\label{d2.11}\ee

\begin{dfn}\label{d2.2}\,\,\,Let $o\in M$ be fixed. The triple $(M, g,o)$ is said to be a {\it \bf model} if for every linear isometry $\psi:$ $M_o\rw M_o$, $\Psi:$
$\exp_o\Sigma(o)\rw\exp_o\Sigma(o)$ is an isometry.
\end{dfn}

 \begin{rem} A conception  of models is introduced in $\cite{GrWu}.$ Definition $\ref{d2.2}$ above is weaker than that in $\cite{GrWu}.$ If $\exp_o:$ $M\rw M$ is a diffeomorphism, they are the same.
\end{rem}

For any $v\in S_o,$ $\gamma(t)=\exp_otv$ is a normal geodesic initiating from the point $o.$ The {\it radial curvature tensor}
along $\gamma(t)$ is a tensor filed of order two, given by
$${\bf R}\Big(\dot{\gamma}(t),X,\dot{\gamma}(t),Y\Big)=\<{\bf R}_{\dot{\gamma}(t)X}\dot{\gamma},Y\>\qfq X,\,\,Y\in M_{\gamma(t)},\quad t\geq0.$$

A model is charactered by its radial curvature. We have

\begin{pro}\label{np2.3} $(M,g,o)$ is a model if and only if there is a function $\kappa$ on $[0,\infty)$ such that
\be {\bf R}\Big(\dot{\gamma}(t),X,\dot{\gamma}(t),Y\Big)=\kappa(t)\<X,Y\>,\quad X,\,\,Y\in M_{\gamma(t)},\label{nn2.11}\ee with $\<X,\dot{\gamma}(t)\>=0$ and $\<Y,\dot{\gamma}(t)\>=0$ for all $t>0,$
where $\gamma(t)=\exp_otv\in \exp_o\Sigma(o),$ for all $v\in M_o$ with $|v|=1.$
\end{pro}

{\bf Proof}\,\,\,Let $(M,g,o)$ be a model and let $\rho>0$ be given. Let $\gamma_i(t)=\exp_otv_i$ with $v_i\in S_0$ for $i=1,$ $2.$
Let $X_i$ be in $M_{\gamma_i(\rho)}$ with $\<X_i,\dot{\gamma}_i(\rho)\>=0$ and $|X_i|=1,$ respectively, for $i=1,$ $2.$
Let $z_i$ be in $S_0$ such that
$$z_i(\rho)=X_i,$$ where $z_i(t)$ are the parallel translations of $z_i$ along $\gamma_i$ such that $z_i(0)=z_i,$ respectively, for $i=1,$ $2.$
Then $$\<z_i,v_i\>=0\qfq i=1,\,2.$$
Suppose $\psi:$ $M_o\rw M_o$ is a linear isometry such that $\psi v_1=v_2$ and $\psi z_1=z_2.$ Then $\Psi,$ given by (\ref{d2.11}), is an isometry  from $\exp_o\Sigma(o)$ to $\exp_o\Sigma(o).$
Then $\Psi(\gamma_1(t))=\gamma_2(t)$ and $\Psi_*z_1(t)=z_2(t).$ We obtain
\beq{\bf R}\Big(\dot{\gamma}_2(\rho),X_2,\dot{\gamma}_2(\rho),X_2\Big)(\gamma_2(\rho))&&={\bf R}\Big(\Psi_*\dot{\gamma}_1(\rho),\Psi_*X_1,\Psi_*\dot{\gamma}_1(\rho),\Psi_*X_1\Big)
(\Psi(\gamma_1(\rho))\nonumber\\
&&={\bf R}\Big(\dot{\gamma}_1(\rho),X_1,\dot{\gamma}_1(\rho),X_1\Big)(\gamma_1(\rho)).\nonumber\eeq Thus, the radial curvatures are the same on the geodesic sphere $\S(\rho),$ centered at $o$ and
with radii $\rho.$ Formula (\ref{nn2.11}) follows with $\kappa$ being the radial curvature.

Conversely, suppose (\ref{nn2.11}) holds. Let $\psi:$ $M_o\rw M_o$ be a linear isometry.
Let $x=\exp_o\rho v\in \exp_o\Sigma(o)$ be given. Let $X_i$ be in $M_x$ with $|X_i|=1$
for $i=1,$ $2.$ Suppose $E_i(t)$ are the parallel translations of some $v_i\in S_o$ along $\gamma(t)=\exp_otv$ such that $E_i(\rho)=X_i,$ respectively, for $i=1,$ $2.$ Let
$$\a_i(t,s)=\exp_ot(v+sv_i)\qfq t\geq0,\,\,s\in\R,\quad i=1,\,\,2.$$ Thus, $J_i(t)=t\exp_{o*}v_i$ are Jacobi fields along $\gamma(t)=\exp_otv.$ By formula (\ref{nn2.11}) and $J'_i(0)=v_i,$ we obtain
\be J_i(t)=f(t)E_i(t)\qfq t\geq0,\,\,i=1,\,\,2,\label{nn2.13}\ee where $f$ is the solution to problem (\ref{nn2.12}) later.
In addition the formulas $\Psi(\a_i(t,s))=\exp_ot\psi(v+sv_i)$ imply
$$\Psi_* J_i(t)=t\exp_{o*}\psi v_i$$ are also Jacobi fields along $\Psi(\gamma(t))=\exp_o t\psi v_i$  for $i=1,$ $2.$ By (\ref{nn2.11}) again,
\be\Psi_*J_i(t)=f(t)\hat{E}_i(t),\label{nn2.14}\ee where $\hat{E}_i(t)$ are parallel translation vector fields along $\Psi(\gamma(t)),$ respectively, for $i=1,$ $2,$ such that
$$\hat{E}_i(0)=\psi v_i\qfq i=1,\,\,2.$$ It follows from (\ref{nn2.13}) and (\ref{nn2.14}) that
$$\<\Psi_*X_1,\Psi_*X_2\>=\<\Psi_*E_1(\rho),\Psi_* E_2(\rho)\>=\<\hat{E}_1(0),\hat{E}_2(0)\>=\<\psi v_1,\psi v_2\>=\<v_1,v_2\>=\<X_1,X_2\>.$$ Thus, $\Psi:$ $\exp_o\Sigma(o)\rw \exp_o\Sigma(o)$ is an isometry.
\hfill$\Box$

\begin{rem} Formula $(\ref{nn2.11})$ means that a model has the same radial curvature on a geodesic sphere $\S(t),$ centered at $o$ with radii
$t>0.$

\end{rem}

Let $(M,g,o)$ be a model and let the radial curvature $\kappa$ be given by (\ref{nn2.11}).  Consider problem
\be \left\{\begin{array}{l} f''(t)+\kappa(t)f(t)=0,\quad t>0,\\
f(0)=0,\quad f'(0)=1.\end{array}\right.\label{nn2.12}\ee
Let $f$ be the solution to problem $(\ref{nn2.12}).$ Then
  \be f(t)=t-\int_0^t(t-s)\kappa(s)f(s)ds,\label{f}\ee which yields
 \be\lim_{t\rw0+}\frac{f(t)}{t}=1.\label{13*}\ee

For our problems here, we need $f$ and $f'$ are all positive functions, for which the following is introduced.
Let the radial curvature $\kappa$ be given in (\ref{nn2.11}). Let
\be\mu_\pm(\lam)=\int_0^\lam s\kappa_\pm(s)ds\qfq \lam>0,\label{kappa13}\ee where
$$\kappa_+(s)=\max\{\kappa(s),\,0\},\quad \kappa_-(s)=\max\{-\kappa(s),\,0\}\qfq s\geq0.$$
We have (\cite{GrWu})

\begin{pro} \label{p1.1}If
\be \mu_+(\lam)\leq1,\label{kappa6}\ee then there exists $0<\mu_0(\lam)\leq1$ such that
\be \mu_0(\lam)\rho\leq f(\rho)\leq e^{\mu_-(\lam)}\rho\qfq\rho\in[0,\lam],\label{kappa4}\ee and
\be \mu_0(\lam)\leq f'(\rho)\leq e^{\mu_-(\lam)}\qfq\rho\in[0,\lam],\label{kappa5*}\ee where $f$ is
the solution to problem $(\ref{nn2.12}).$
\end{pro}

{\bf Proof}\,\,\,Let $p_+$ and $p_-$ solve the problems
\be\left\{\begin{array}{l}p_+''+\kappa_+ p_+=0\qfq \rho>0,\\
p_+(0)=0,\quad p_+'(0)=0\end{array}\right.\label{kappa1}\ee and
\be\left\{\begin{array}{l}p_-''-\kappa_- p_-=0\qfq \rho>0,\\
p_-(0)=0,\quad p_-'(0)=0,\end{array}\right.\label{kappa2}\ee respectively.

Let $\eta(\rho)=p_-/p_-'$ for $\rho>0.$ By (\ref{kappa2}) we have
$$\eta'=1-\kappa_-\eta\leq1,\quad\mbox{and, then}\quad \eta\leq\rho\qfq\rho\geq0,$$ since $\eta(0)=0.$ It follows that
$$\frac{p_-''}{p_-'}=\kappa_-\eta\leq\kappa_-(\rho)\rho\qfq\rho\geq0.$$
Integrating the above inequality over $(0,\rho)$ yields
\be 1\leq p_-'(\rho)\leq e^{\mu_-(\lam)}\qfq\rho\in[0,\lam],\label{kappa3}\ee which implies
\be \rho\leq p_-(\rho)\leq e^{\mu_-(\lam)}\rho\qfq\rho\in[0,\lam].\label{kappa5}\ee

On the other hand, from (\ref{kappa1}) we obtain
\be  \mu_0(\lam)\leq p_+'(\rho)=1-\int_0^\rho\kappa_+(s)p_+(s)ds\leq1,\label{kappa10}\ee
and then \be\mu_0(\lam)\rho\leq p_+(\rho)\leq\rho\qfq\rho\in[0,\lam],\label{kappa7}\ee where $\mu_0(\lam)=\min_{0\leq\rho\leq\lam}p'_+(\lam)\leq1.$ Moreover, we claim that  assumption
(\ref{kappa6}) implies that $\mu_0(\lam)>0.$ Otherwise, if there were a point $\rho_0\in[0,\lam]$ such that $p_+'(\rho_0)=0,$ then by (\ref{kappa7})
$$1=\int_0^{\rho_0}\kappa_+(s)p_+(s)ds<\mu_+(\lam)\leq1,$$ a contradiction.

By a comparation argument for ordinary differential equations we obtain
\be\frac{p_+'}{p_+}\leq\frac{f'}{f}\leq\frac{p_-'}{p_-}\qfq\rho\in[0,\lam].\label{kappa8}\ee
Integrating the above inequalities over $(\varepsilon,\rho)$ gives
$$p_+\frac{f(\varepsilon)}{p_+(\varepsilon)}\leq f\leq \frac{f(\varepsilon)}{p_-(\varepsilon)}p_-\qfq0<\varepsilon\leq\rho\leq\lam.$$
Letting $\varepsilon\rw0+$ in the above inequalities we have
\be p_+\leq f\leq p_-\qfq\rho\in[0,\lam],\label{kappa9}\ee since $$\lim_{\rho\rw0+}\frac{f(\rho)}{\rho}=\lim_{\rho\rw0+}\frac{p_\pm(\rho)}{\rho}=1.$$
Then (\ref{kappa4}) follows from (\ref{kappa9}), (\ref{kappa7}), and (\ref{kappa5}).

Finally, (\ref{kappa5*}) follows from (\ref{kappa8}), (\ref{kappa9}), (\ref{kappa3}), and (\ref{kappa10}).\hfill$\Box$\\

The following result  is immediate from Proposition \ref{p1.1}  which is an improvement of Lemmas 4.5 and 4.6 in \cite{GrWu}.

\begin{pro} \label{np1.1}If
\be \mu_+(\infty)\leq1,\quad \mu_-(\infty)<\infty,\label{kappa11}\ee then there are $0<\mu_0\leq\mu_1<\infty$ such that
\be \mu_0\leq f'\leq \mu_1\quad\mbox{and}\quad \mu_0\rho\leq f\leq \mu_1\rho\qfq\rho\in[0,\infty).\label{pappa12}\ee
\end{pro}

Let $\rho=\rho(x)$ be the distance function in the metric $g$ from $o$ to $x\in M.$
We recollect some properties of a model from \cite{GrWu} in the following.

\begin{pro}\label{np2.4} Let $(M,g, o)$ be a model. Then

$(i)$\,\,\,Any Jacobi field $J(t)$ is in the form
$$J(t)=f(t)E(t)\qfq t>0,$$ where $E(t)$ is the parallel translation.

$(ii)$\,\,\,The Hessian of the distance function $\rho$ is given by
$$D^2\rho=\frac{f'(\rho)}{f(\rho)}(g-D\rho\otimes D\rho)\qfq \rho(x)>0,\,\,x\in\Sigma(o).$$

$(iii)$\,\,\, In the geodesic polar coordinates the metric $g$  has the expression:
 \be g=d\rho^2+f^2(\rho)d\theta^2\qfq \rho>0,\,\,x\in \Sigma(o).\label{12*}\ee
\end{pro}

 \begin{pro} \label{np2.5} Let $(M,g,o)$ be a model. Let $\u$ be a radical deformation
 given by $(\ref{7}).$
  Then \be \Big(d\u\Big)=\diag\Big(\var'(\rho),\tau(\rho),\cdots,\tau(\rho)\Big)\qfq x\in \Sigma(o),\label{13}\ee where
\be\tau(\rho)=\frac{f\circ\var(\rho)}{f(\rho)}\qfq x\in\Sigma(o).
\label{14}\ee
\end{pro}

{\bf Proof}\,\,\,Let $x=\exp\rho v$ where $v\in M_o$ with $|v|=1.$ Let $\{E_i(t)\}$ be the parallel translation orthonormal basis of $M_{\gamma(t)}$ along $\gamma(t)=\exp_otv$ with $E_1=v.$ Then
$$ E_1(t)=D\rho(\gamma(t))\qfq t>0.$$
By Proposition \ref{pl3} and Proposition \ref{np2.4} (i), we have
$$\<E_1(\var(\rho)),\u_*E_1(\rho)\>=\<D\rho(\var(\rho)),\u_*D\rho(x)\>=\var'(\rho),$$
$$ \<E_1(\var(\rho)),\u_*E_j(\rho)\>=0\qfq 2\leq j\leq n,$$
$$\<E_i(\var(\rho)),\u_*E_j(\rho)\>=\<E_i(\var(\rho)),\frac{1}{f(\rho)}\u_*J_j(\rho)\>=\frac{f\circ\var(\rho)}{f(\rho)}\delta_{ij},$$ for
$2\leq i,\,\,j\leq n.$ \hfill$\Box$

\subsection{Compressible Case}
\hskip\parindent Let $\u:$ $M\rw M$ be a differentiable map. In order to compute equilibrium equations to  energy (\ref{4}), we consider
the vector bundle $\zeta=\u^{-1}TM$ over  the base manifold
$(M, g)$, induced by the map $\u$,
\be\zeta=\cup_{x\in M}M_{\u(x)}.\label{18**}\ee  The projection map
$\pi:\zeta\rw M$ is given by $$\pi(x,Y)=x\qfq x\in M,\quad
Y\in M_{\u(x)}.$$ In local coordinates a {\it section} $H$ of $\zeta$ is in the form \be
H(x)=\sum_{i=1}^nh_i(x)\pl_{x_i}|_{\u(x)}\qfq x\in M,\label{15}\ee
where $h_i\in C^\infty(M)$ for all $i$. Denote by  $\Ga(\zeta)$
all sections of $\zeta.$ The connection $D:
\X(M)\times\Ga(\zeta)\rw\Ga(\zeta)$, induced by the metric $g$,
is given by \be
D_XH=\sum_{i=1}^n[X(h_i)\pl_{x_i}|_{\u(x)}+h_i(x)(D_{\u_*X}\pl_{x_i})\circ\u],
\label{16}\ee where $X\in\X(M)$ and $H\in\Ga(\zeta)$. Furthermore, for $H_1,$ $H_2\in\zeta$ and $X\in \X(M),$ we have
$$ X\<H_1,H_2\>=\<D_{X}H_1,H_2\>+\<H_1,D_{X}H_2\>.$$

For $H\in\Ga(\zeta)$, we set
\be{\bf v}(t)(x)=\exp_{\u(x)}tH(x)\qfq t\geq0,\quad
x\in M,\label{17}\ee where $\exp_{\u(x)}$: $M_{\u(x)}\rw
M$ is the exponential map in the metric $g$ at $\u(x)$ along
the vector $H(x)\in M_{\u(x)}$ for each $x\in M$. For each $x$
fixed, ${\bf v}(t)$ is a geodesic on $M$ initiating from $\u(x)$
in the metric $g$ and
$${\bf v}(0)(x)=\u(x),\quad \dot{\bf v}(0)(x)=H(x)\in M_{\u(x)}\qfq x\in M.$$
${\bf v}(t)$ is said to be a {\it variation} of $\u$ (Chapter $8$ of
\cite{Jo}).

Let $x\in M$ be given.  For $e\in
M_x$, ${\bf v}_*(t)e$ is a vector field along the geodesic ${\bf
v}(t)(x)$.  We have

\begin{lem}\label{l5} Let $\u:$ $M\rw M$ be a differentiable map
and let ${\bf v}(t)$ be a variation of $\u$, given by $(\ref{17})$
for $t\in(-\varepsilon,\varepsilon)$. Then \be D_{\dot{{\bf
v}}(0)}{\bf v}_*e=D_eH\qfq e\in M_x,\quad
x\in M.\label{18}\ee
\end{lem}

{\bf Proof}\,\,\,We do a computation in local coordinates
$x=(x_1,\cdots,x_n)$. Let $H\in\Ga(\zeta)$ be given by (\ref{15}). Let
$e=\sum_{i=1}^n\a_i\pl_{x_i}|_x$ and
$${\bf v}(t)(x)=(v_1(t,x),\cdots,v_n(t,x))\qfq(t,x)\in(-\varepsilon,\varepsilon)\times M,$$
where $$(v_1(0,x),\cdots,v_n(0,x))=\u(x),\quad \dot{\v}(0)x=H(x)\qfq x\in M.$$ Then
$$\dot{\bf v}(t)=\sum_{i=1}^n \dot{v}_{i}(t)\pl_{x_i}|_{{\bf v}(t)},
\quad{\bf v}_*(t)\pl_{x_j}=\sum_{i=1}^nv_{ix_j}\pl_{x_i}|_{{\bf
v}(t)},$$ where $\dot{v}_{i}(0)=h_i(x)$ and $\v_*(0)\pl x_i=\u_*\pl x_i$ for $1\leq i\leq n.$ We have
\beq D_{\dot{{\bf v}}(0)}{\bf v}_*e&&=D_{\dot{{\bf
v}}(0)}[\sum_i(\sum_j\a_jv_{ix_j})\pl_{x_i}|_{{\bf
v}(t)}]\nonumber\\
&&=\sum_i[\sum_j\dot{\bf
v}(0)(\a_jv_{ix_j})\pl_{x_i}|_{\u(x)}+\sum_j\a_jv_{ix_j}\sum_kh_k(D_{\pl_{x_k}}\pl_{x_i})|_{\u(x)}]\nonumber\\
&&=\sum_l[\sum_j\a_j(v_{lx_jt}(0)+\sum_{ik}v_{ix_j}(0)h_k\Ga_{ki}^l\circ\u)]\pl_{x_l}|_{\u(x)}\nonumber\\
&&=\sum_l[\sum_j\a_j(h_{lx_j}+\sum_{ik}u_{ix_j}h_k\Ga_{ki}^l\circ\u)]\pl_{x_l}|_{\u(x)}.\label{19}\eeq
It follows  from (\ref{19}) that \beq
D_eH&&=\sum_i[e(h_i)\pl_{x_i}|_{\u(x)}+h_iD_{\u_*e}\pl_{x_i}]\nonumber\\
&&=\sum_{ij}\a_j[h_{ix_j}\pl_{x_i}|_{\u(x)}+h_i\sum_ku_{kx_j}(D_{\pl_{x_k}}\pl_{x_i})|_{\u(x)}]\nonumber\\
&&=D_{\dot{{\bf v}}(0)}{\bf v}_*e.\nonumber\eeq \hfill$\Box$\\

Let $e_1,$ $\cdots,$ $e_n$ be an orthonormal basis of $M_o$. Consider the usual pole coordinates $(\rho,\theta)$ on the
Euclidean space $\R^n$, given by
\be\left\{\begin{array}{l} z_1=\rho\cos\theta_1,\\
z_2=\rho\sin\theta_1\cos\theta_2\\ \cdots\cdots\\
z_{n-1}=\rho\sin\theta_1\cdots\sin\theta_{n-2} \cos\theta_{n-1}\\
z_n= \rho\sin\theta_1\cdots\sin\theta_{n-2} \sin\theta_{n-1},
\end{array} \right.\label{gg1}\ee where $\theta=(\theta_1,\cdots,\theta_{n-1})$
and
$$0\leq\theta_1\leq\pi,\quad\cdots,\quad0\leq\theta_{n-2}\leq\pi,\quad
0\leq\theta_{n-1}\leq2\pi.$$ Let
\be\pi(\theta)=\frac{1}{\rho}\sum_{i=1}^nz_ie_i.\label{gg4}\ee We have
$$\lim_{\rho\rw0+}\int_{S_o}\<X,D\rho\>d\theta=\lim_{\rho\rw0+}\int_{S_o}\<X,\exp_{o*}\pi(\theta)\>d\theta=\int_{S_o}\<X(o),\pi(\theta)\>d\theta=0,$$ and, thus,
\beq\lim_{\rho\rw0}\int_{S_o}\frac{1}{\rho}\<X,D\rho\>d\theta&&=\lim_{\rho\rw0+}\int_{S_o}\frac{\<X,D\rho\>-\<X(o),\pi(\theta)\>(o)}{\rho}d\theta\nonumber\\
&&=\int_{S_o}\<D_{\pi(\theta)}X,
\pi(\theta)\>d\theta.\label{nnn2.42}\eeq

Let the function $W:M_+^{n\times n}\rw\R$ satisfy assumption
(\ref{2}). It is well known (\cite{TrNo}) that  there exists a
symmetric function $$\Phi:\R_{++}^n\rw\R,\quad\R^n_{++}=\{\,
c=(c_1,\cdots,c_n)\in\R^n,\,\,c_i>0\,\,\mbox{for}\,\,1\leq i\leq
n\,\},$$ such that\be W(F)=\Phi(v_1,\cdots,v_n)\quad\mbox{for
all}\quad F\in M_+^{n\times n},\label{2.15}\ee where $v_1,$
$\cdots,$ $v_n$ denote the singular values (or principal stretches)
of $F$ (i.e., the eigenvalues of $(F^\tau F)^{1/2}$). It is known
(\cite{Ball}) that $W\in C^r(M_+^{n\times n})$ if and only if
$\Phi\in C^r(\R^n_{++})$ for $r=0,$ $1,$ $2$ or $\infty$. We write
$\Phi_i=\dfrac{\pl\Phi}{\pl v_i},$ etc. If $W\in C^1(M_+^{n\times
n})$ and if $F=\diag(v_1, \cdots, v_n),$ $v_i>0$, then
\be\frac{\pl W}{\pl F}(F)=\diag(\Phi_1,\cdots,\Phi_n),\label{n2.26}\ee where $\Phi_i=\Phi_i(v_1,\cdots,
v_n)$ for all $i$. Moreover, the symmetry of the function $\Phi$
implies
\be\Phi_i(v_1,v,\cdots,v)=\Phi_2(v_1,v,\cdots,v)\qfq 2\leq i\leq n.\label{n2.28}\ee

Let $\u:$ $\Om\rw M$ be a deformation. We say that $\u\in W^{1,p}(\Om,M)$ if
$$\int_\Om |F|^pdg<\infty,$$ where $F$ is given by (\ref{4*}) and $1\leq p<\infty.$ Moreover, we define
$$\det d\u(x)=\det F\qfq x\in\Om.$$

Let $\Om=\B$ be the unit geodesic ball centered at $o.$ Let $H\in\Ga(\zeta)$ be a section. Then the differential of $H$ in the connection $D$ can be
defined by
$$DH(x)=\Big(\<E_i,D_{e_j}H\>\Big)\qfq x\in M,$$ where $\{e_i\}$ and $\{E_i\}$ are orthonormal bases of $M_x$ and $M_{\u},$ respectively.
We say that a deformation $\u\in W^{1,1}(\B,M)$ is {\it an equilibrium solution}
if $\det d\u>0,$ a.e. $x\in\B,$
$\frac{\pl W}{\pl F_{ij}}\in L^1(\B)$ for $1\leq i,\,j\leq n,$ and
$$\int_\B \<D_FW,DH\>dg=0\quad\mbox{for all}\quad H\in C_0^\infty(\B,\Ga(\zeta)).$$

\begin{rem} In general the vector bundle $\Ga(\zeta)$, given by $(\ref{18**}),$ may depend on the deformation $\u.$ Then does $DH.$ If $M=\R^n$ with the Euclidean metric,
then $M_{\u(x)}=M_x=\R^n$ and a section $H$ and its differential $DH$  are independent of the deformation.
\end{rem}

\begin{thm}\label{nt1.1} Let $(M,g,o)$ be a model  with $\mu_+(1)\leq1.$ Let $\u$ be a radical deformation
 given by $(\ref{7})$ such that $\mu_+(\var(1))\leq1.$  Then $\u$ is an
equilibrium solution to problem $(\ref{5})$ if and only if $\var\in W^{1,1}(0,1),$
$\var'(\rho)>0$ a.e. $\rho\in(0,1],$ $f^{n-1}\Phi_1,$ $f^{n-1}\Phi_2\in L^1(0,1),$ and\be
[f^{n-1}(\rho)\Phi_1(\rho)]_\rho=(n-1)f^{n-2}(\rho)f'\circ\var(\rho)\Phi_2(\rho)\qfq x\in\Om,\quad
\rho(x)>0,\label{22*}\ee where
$$\Phi_i(\rho)=\Phi_i(\var',\tau,\cdots,\tau)\qfq \rho>0,\quad i=1,\,\,2,$$
and $\tau$ is given by $(\ref{14}).$
\end{thm}

{\bf Proof}\,\,\,Let $H$ be a section of $\Ga(\zeta)$ with a compact
support such that $\supp H\subset\B$ and let ${\bf v}(t)=\exp_{\u(x)}tH(x)$ be a variation of $\u$. Let
$x=\exp_o\rho v\in \Sigma(o)$ be given, where $v\in M_o$ with $|v|=1.$  Let $\{E_i\}$  be   an orthonormal basis of $M_o$ with $E_1=v.$ We transport  $\{E_i\}$ along
the geodesic $\gamma(t)=\exp_otv$ parallelly to obtain the orthonormal bases $\{E_i(t)\}$ of $M_{\gamma(t)}$ for $t\geq0.$
Next, we  transport parallelly  the orthonormal basis $\{E_i(\var(\rho))\}$ of $M_{\u(x)}$ along the geodesic $\v(t)=\exp_{\u(x)}tH$  to have the orthonormal bases $\{\hat{E}_i(t)\}$
of $M_{\v(t)}$ for which
\be \hat{E}_i(0)=E_i(\var(\rho))\qfq 1\leq i\leq n.\label{22}\ee
In particular, \be \hat{E}_1(0)=E_1(\var(\rho))=D\rho(\var(\rho)).\label{25*}\ee
By  Propositions  \ref{pl3} and Proposition \ref{np2.4} (i), we have
\be\u_*E_1(\rho)=\var'(\rho)D\rho(\var(\rho)),\quad \u_*E_i(\rho)=\tau(\rho) E_i(\var(\rho))\qfq 2\leq i\leq n.\label{23}\ee

 Let
$$F(t)=\Big(\<\hat{E}_i(t),{\bf v}_*(t)E_j(\rho)\>\Big).$$
By Proposition \ref{np2.5}, $$F(0)=\diag\Big(\var',\tau,\cdots,\tau\Big),$$ since
$\v_*(0)=\u_*.$ We obtain, by  Lemma
\ref{l5}, (\ref{22}), (\ref{n2.26}), (\ref{n2.28}) and (\ref{25*}), \beq \frac{d}{dt}W(d{\bf
v})\Big|_{t=0}&&=\sum_{ij=1}^n\frac{\pl W}{\pl
F_{ij}}(F(0))\frac{\pl\<\hat{E}_i(t),{\bf v}_*(t)E_j(\rho)\>}{\pl
t}\Big|_{t=0}\nonumber\\
&&=\sum_{ij=1}^n\frac{\pl W}{\pl
F_{ij}}(F(0))\<\hat{E}_i(0),D_{E_j(\rho)}H\>\circ\u(x)\nonumber\\
&&=\Phi_1(\rho)\<D\rho,D_{D\rho}H\>\circ\u(x)+\sum_{i=2}^n\Phi_i
\<E_i(\var(\rho)),D_{E_i(\rho)}H\>\circ\u(x)\nonumber\\
&&=\frac{\pl}{\pl\rho}[\Phi_1(\rho)\<H,D\rho\>]-[\frac{\pl}{\pl\rho}\Phi_1(\rho)]
\<H,D\rho\>\nonumber\\\
&&\quad+\Phi_2(\rho)\sum_{i=2}^n
\<E_i(\var(\rho)),D_{E_i(\rho)}H\>\circ\u(x).\label{24}\eeq

Next, let \be H=h_0(x)D\rho|_{\u(x)}+\tilde{H},\label{nn2.36}\ee where
$$\tilde{H}=\sum_{j=1}^{n-1}h_j(x)\pl\theta_j\Big|_{\u(x)}$$ and $(\rho,\theta)$ are the geodesic polar coordinates on $M$ initiating from $o$ in the metric $g.$
It follows from (\ref{nn2.36}) that $\<H,D\rho\>=h_0(x)$ for $x\in M,$ i.e., $\<H,D\rho\>$ is a function on $M,$ which is independent of deformation $\u.$

Using
(\ref{23}) and Proposition \ref{np2.4} (ii), we obtain \beq
\<E_i(\var(\rho)),D_{E_i(\rho)}H\>&&=\<E_i(\var(\rho)),\quad E_i(\rho)(h_0)D\rho|_{\u(x)}+
h_0(x)D_{\u_*E_i(\rho)}D\rho+D_{E_i(\rho)}\tilde{H}\>\nonumber\\
&&=\<H,D\rho\>\tau(\rho)D^2\rho(E_i(\var(\rho)),E_i(\var(\rho)))+\<E_i(\var(\rho)),\tilde{D}_{E_i(\rho)}\tilde{H}\>\nonumber\\
&&=\frac{\<H,D\rho\>f'\circ\var(\rho)}{ f(\rho)}+\<E_i(\var(\rho)),\tilde{D}_{E_i(\rho)}\tilde{H}\>,
\label{25}\eeq for $2\leq i\leq n,$ where $\tilde{D}$ is the induced connection of
$\S(\var(\rho))$ from the metric $g$ and where $\S(\var(\rho))$ denotes the geodesic sphere with radii $\var(\rho)$ centered at $o.$

Let functions $p_j$ on $\S(\var(\rho))$ be defined by
$$p_j(\u(x))=h_j(x)\qfq\u(x)\in\S(\var(\rho)),\quad 1\leq j\leq n-1.$$
Using  (\ref{16}), we obtain
\beq\tilde{D}_{E_i(\rho)}\tilde{H}\Big|_{\u(x)}&&=\sum_{j=1}^{n-1}[E_i(h_j)(x)\pl\theta_j\Big|_{\u(x)}+h_j(x)\tilde{D}_{\u_*E_i(\rho)}\pl\theta_j\circ\u(x)]\nonumber\\
&&=\tau\sum_{j=1}^n[E_i(p_j)(\u)\pl\theta_j\Big|_{\u}+p_j(\u)\tilde{D}_{E_j(\var(\rho))}\pl\theta_j\circ\u].\nonumber\eeq Thus,
\be\sum_{i=2}^{n}\<E_i(\var(\rho)),D_{E_i(\rho)}\tilde{H}\>\circ\u(x)=\tau\tilde{\div}\tilde{H}_0\circ \u(x),\label{n2.33}\ee where $\tilde{\div}$ is the divergence of the induced metric on
$\S(\var(\rho))$ from the metric $g$ and $\tilde{H}_0$ is a vector field on $\S(\var(\rho))$, given by
\be\tilde{H}_0=\sum_{j=1}^{n-1}p_i(\u)\pl\theta_j\Big|_{\u(x)}\qfq \u(x)\in \S(\var(\rho)).\label{nn2.39}\ee

 By inserting (\ref{25}) and (\ref{n2.33}) into
(\ref{24}), we have \beq \frac{d}{dt}W(d{\bf v})\Big|_{t=0}&&
=\frac{\pl}{\pl\rho}[\Phi_1(\rho)\<H,D\rho\>]+[(n-1)\frac{f'\circ\var(\rho)}{f(\rho)}\Phi_2(\rho)-
\frac{\pl}{\pl\rho}\Phi_1(\rho)]
\<H,D\rho\>\nonumber\\
&&\quad+\tau(\rho)\Phi_2(\rho)(\tilde{\div}\tilde{H}_0)\circ\u(x).\label{26}\eeq

Using (\ref{12*}) and (\ref{26}), we obtain \beq
&&\int_\B\<D_FW,DH\>dg=\int_\B\frac{d}{dt}W(d{\bf
v})\Big|_{t=0}dg=\int_0^1\int_{S_o}\frac{d}{dt}W(d{\bf
v})\Big|_{t=0}f^{n-1}(\rho)d\rho
d\theta\nonumber\\
&&=\int_0^1\Big\{(n-1)f^{n-2}(\rho)f'\circ\var(\rho)\Phi_2(\rho)
-[f^{n-1}(\rho)\Phi_1(\rho)]_\rho\Big\}d\rho
\int_{S_o}\<H,D\rho\>d\theta\nonumber\\
&&\quad+\int_0^1\tau^{2-n}
\Phi_2(\rho)d\rho\int_{\S(\var(\rho))}\tilde{\div}\tilde{H}_0d\tilde{g}
\nonumber\\
&&\quad+\int_{S_o}(f^{n-1}\Phi_1\<H,D\rho\>d\theta\Big|_{\rho=1}
-\lim_{\varepsilon\rw0}f^{n-1}(\varepsilon)\Phi_1(\varepsilon)\int_{S_o}\<H,D\rho\>d\theta,\label{28}\eeq
where $\tilde{g}$ is the induced metric on $\S(\var(\rho))$ from the
metric $g$ and $S_o\subset M_o$ is the unite sphere of $M_o.$

Since $(\S(\var(\rho)),\tilde{g})$ is a compact manifold without a
boundary, the second integral in the right hand side of (\ref{28})
is zero.

Let radical deformation $\u$ be an  equilibrium solution.
We take $H\in\Ga(\zeta)$ with $\supp H\subset\subset \B/\{o\}.$  Then the last two integrals in the right hand (\ref{28}) are zero.
Thus, equation (\ref{22*}) follows from (\ref{28}).

Conversely, suppose that $f^{n-1}\Phi_1,$ $f^{n-1}\Phi_2\in L^1(0,1)$ and that  equation (\ref{22*}) is true. It follows from  equation (\ref{22*}) and
the relation (\ref{28}) that for any $H\in C_0^\infty(\B,\Ga(\zeta))$
\be\int_\B\<D_FW,DH\>dg=-\lim_{\varepsilon\rw0}f^{n}(\varepsilon)\Phi_1(\varepsilon)\frac{\varepsilon}{f(\varepsilon)}\int_{S_o}\frac{\<H,D\rho\>}{\varepsilon}d\theta.\label{36}\ee
Since $f^{n-1}\Phi_1,$ $f^{n-1}\Phi_2\in L^1(\B),$ it follows from  (\ref{kappa5*}) and  (\ref{22*}) that
$$(f^n\Phi_1)_\rho=f'f^{n-1}\Phi_1+(n-1)f'\circ\var f^{n-1}\Phi_2\in L^1(0,1).$$ Thus, $\lim_{\rho\rw0+}f^n\Phi_1=0.$
By (\ref{36}) and (\ref{nnn2.42}), $\u$ is an equilibrium solution.  \hfill$\Box$

\begin{rem}
The above theorem is Theorem $4.2$ in $\cite{Ball1}$ for the
isotropic materials if $M=\R^n$ with the Euclidean metric. In that case $f(\rho)=\rho$ and
$$\Big(d\u\Big)=\Big(\var'(\rho),\frac{\var(\rho)}{\rho},\cdots,\frac{\var(\rho)}{\rho}\Big).$$
Equation $(\ref{22*})$ becomes
$$[\rho^{n-1}\Phi_1(\rho)]_\rho=(n-1)\rho^{n-2}\Phi_2(\rho).$$
\end{rem}

Throughout this paper a function $\var\in W^{1,1}(0,1)$ is said to {\it an equilibrium solution} to  problem (\ref{22*}) if $\var$ satisfies equation (\ref{22*}) and $\var'(\rho)>0$
for $\rho\in(0,1],$ and is such that $f^{n-1}\Phi_1,$ $f^{n-1}\Phi_2\in L^1(0,1).$ By a similar argument in \cite{Ball1}, we have

\begin{thm}\label{t1.2} Let $$\Phi_{11}(v_1,v_2,\cdots,v_2)>0\qfq v_1>0,\,\,v_2>0.$$
If $\var\in W^{1,1}(0,1)$ is an equilibrium solution to  problem $(\ref{22*}),$ then $\var\in C^1(0,1].$
\end{thm}

\subsection{Incompressible Case}
\hskip\parindent
A map $\u:$ $M\rw M$ is said to be {\it incompressible} if, for any
element of volume $\omega$ of $M_{\u(x)}$, $\u^*\omega$ is an element of volume of $M_x$ for $x\in M.$ Let $(M,g,o)$ be a model and let $\u$ be a radical deformation
 given by $(\ref{7}).$ Let $x=\exp_o\rho v\in\Sigma(o)$ with $v\in M_o$ and $|v|=1.$ Let $\{e_i\}$ be an orthonormal basis of $M_x$ with the positive orientation such that $e_1=D\rho(x).$ Let $\omega$ be a volume element of $M_{\u(x)}.$ By Proposition \ref{np2.4} (i), we have
 $$\u^*\omega(e_1,e_2,\cdots,e_n)=\omega(\u_*e_1,\u_*e_2,\cdots,\u_*e_n)
 =\var'(\rho)\tau^{n-1}(\rho),$$ where $\tau$ is given by (\ref{14}). Thus, $\u$ is incompressible if and only if
 $$\var'(\rho)f^{n-1}\circ\var(\rho)=f^{n-1}(\rho)\qfq \rho>0.$$
We define
\be\si(t)=\int_0^tf^{n-1}(s)ds\qfq t\geq0.\label{si}\ee
 Then only possible such deformations satisfy
\be \si(\var(\rho))=\int_0^\rho f^{n-1}(s)ds+\eta,\label{31}\ee where $\eta$ is a constant. To get $\var$ from the above equation, we need some assumptions on the radical curvature.

Let conditions  (\ref{kappa11}) hold true. It follows from (\ref{31}) that for an incompressible deformation
 \be\var(\rho)=\si^{-1}\Big(\si(\rho)+\si(A)\Big),\label{39*}\ee where $A=\var(0).$

 \begin{lem} Let $(M,g,o)$ be a model  with $\mu_+(1)\leq1.$
 Then  the radical deformations with $\mu_+(\var(1))\leq1,$ which are incompressible, belong to $W^{1,p}(B,M)$ for $1\leq p<n.$
 \end{lem}

 {\bf Proof}\,\,\,We need to prove
 \be \int_0^1\Big[\var'^2(\rho)+(n-1)\tau^2(\rho)\Big]^{p/2}f^{n-1}(\rho)d\rho<\infty\qfq 1\leq p<n.\label{32}\ee
 If $\var(0)=0,$ then $\var(\rho)=\rho.$ The above estimate is trivial. Let us assume that $\var(0)=A>0.$ Then estimate (\ref{32}) follows from (\ref{13*}). \hfill$\Box$\\

Let $\u:$ $M\rw M$ be a deformation. The determinant of $d\u$ is given by
$$\det d\u(x)=\det\Big(\<E_i,\u_*e_i\>\Big)\qfq x\in M,$$ where $\{e_i\}$ and $\{E_i\}$ are orthonormal bases of $M_x$ and $M_{\u(x)}$, respectively.

Let $W$ be a constitutive function satisfying (\ref{2}) and let $\B$ be the unit geodesic ball centered at $o.$
The equilibrium equations for incompressible radial deformations are the Euler-Lagrange equations for the functional
$$I(\u)=\int_\B\{W(d\u)-p(x)[\det d\u-1]\}dg,$$ where the pressure $p(x)$ is a Lagrange multiplier corresponding to the constraint of incompressibility. Then
a deformation $\u:$ $\B\rw M$ is said to be an equilibrium solution (\cite{Ball1}) with corresponding measurable pressure $p(x)$ if $\det\u=1$ a. e. in $\B,$
$$\pl W(d\u)/\pl F_{ij}-p(x)(\adj d\u)_{ij}\in L^1(\B)\qfq1\leq i,\,\,j\leq n,$$ and
$$\int_\B\<D_FW-pD_F\det d\u,\,\,DH\>dg=0\qfq H\in C_0^\infty(\B,\Ga(\zeta)).$$

Let
$$\hat{\Phi}(v)=\Phi(v^{1-n},v,\cdots,v)\qfq v>0.$$ The following theorem is Theorem 4.3 in \cite{Ball1} if $M=\R^n$ is the Euclidean space.

\begin{thm} \label{t2}Let $(M,g,o)$ be a model  with $\mu_+(1)\leq1.$ The radical deformation $(\ref{39*})$ with $A>0$ and $\mu_+(\var(1))\leq1$ is an equilibrium solution if and only if
\be\frac{\tau^{n-1}}{(\tau^n-1)^2}\hat{\Phi}'(\tau)\in L^1(\delta,\infty)\qfq \delta>1.\label{42} \ee
 In this case the corresponding pressure is given by
\be p=\int^{1}_{\rho}
\frac{f'\circ\var(\rho)}{f\circ\var(\rho)} \tau^{2-n}(\rho)\hat{\Phi}'(\tau(\rho))d\rho+\tau^{1-n}(\rho)\Phi_1(\rho)
+c,\label{43*}\ee
where $\Phi_1(\rho)=\Phi_1(\tau^{1-n},\tau,\cdots,\tau),$ $\tau=\tau(\rho)$ is given by $(\ref{14}),$ and  $c$ is a constant.
\end{thm}

{\bf Proof}\,\,\,Let $p$ be defined by (\ref{43*}). It is easy to check that $X=\Phi_1-p\tau^{n-1}$ satisfies the equation
\be X_\rho+(n-1)(\frac{f'}{f}-\frac{\var'f'\circ\var}{f\circ\var})X-\frac{f'\circ\var}{f}\hat{\Phi}'(\tau)=0\qfq \rho>0,\label{44n}\ee
where $p$ is given by (\ref{43*}).

Let $H$ be a section of $\Ga(\zeta)$ with a compact
support on $\B$ and let ${\bf v}(t)$ be a variation of $\u$, given by (\ref{17}). Let
$x=\exp_o\rho v\in \Sigma(o)$ be given, where $v\in M_o$ with $|v|=1.$  Let $\{E_i\}$  be   an orthonormal basis of $M_o$
with $E_1=v.$  We transport  $\{E_i\}$ along
the geodesic $\gamma(t)=\exp_otv$ parallelled to obtain the orthonormal bases $\{E_i(t)\}$ of $M_{\gamma(t)}$ for $t\geq0.$
 Next, we  transport parallelly  the orthonormal basis $\{E_i(\var(\rho))\}$ of $M_{\u(x)}$ along the geodesic $\v(t)=\exp_{\u(x)}tH$
 to have the orthonormal bases $\{\hat{E}_i(t)\}$
of $M_{\v(t)}$ for $t\geq0$ such that the relations (\ref{22}) , (\ref{25*}) and (\ref{23}) hold.

Denote
$$ P_i(t)=\Big(\<\hat{E}_i(t),\v_*(t)E_1\>,\cdots,\<\hat{E}_i(t),\v_*(t)E_n\>\Big)^T\qfq1\leq i\leq n.$$ By (\ref{22}), (\ref{18}),
 (\ref{23}), (\ref{nn2.36}), (\ref{25}) and (\ref{n2.33}), we have
\beq\frac{\pl P_i}{\pl t}\Big|_{t=0}&&=\Big(\<\hat{E}_i(0),D_{\v(0)}\v_*E_1\>,\cdots,\<\hat{E}_i(0),D_{\v(0)}\v_*E_n\>\Big)^T
\nonumber\\
&&=\Big(\<E_i\circ\var(\rho),D_{E_1}H\>,\cdots,\<E_i\circ\var(\rho),D_{E_n}H\>\Big)^T\qfq1\leq i\leq n\nonumber\eeq
and
\beq \frac{\pl\det d\v}{\pl t}\Big|_{t=0}&&=\sum_{k=1}^n\det\Big(P_1(0),\cdots,\dot{P}_k(0),\cdots, P_n(0)\Big)\nonumber\\
&&=\tau^{n-1}\<D\rho,D_{D\rho}H\>(\u(x))
+\tau^{n-2}\var'\sum_{i=2}^n\<E_i\circ\var,D_{E_i(\rho)}H\>\nonumber\\
&&=\tau^{n-1}\<D\rho,D_{D\rho}H\>(\u(x))+(n-1)\frac{f'\circ\var}{f\circ\var}\<H,D\rho\>\nonumber\\
&&\quad+\tilde{\div}\tilde{H}_0\Big|_{\u(x)},\quad\label{39}\eeq where $\tilde{H}_0$ is a vector field on $\S(\var),$ given by (\ref{nn2.39}).
From (\ref{26}), (\ref{39}) and (\ref{44n}), a similar argument as in the proof of Theorem \ref{nt1.1} yields
\beq \frac{\pl I(\v)}{\pl t}\Big|_{t=0}&&=\int_{S_o}f^{n-1}(\Phi_1-p\tau^{n-1})\<H,D\rho\>d\theta\Big|_0^1\nonumber\\
&&\quad-\int_0^1f^{n-1}\{X_{\rho}+(n-1)(\frac{f'}{f}-\frac{\var'f'\circ\var}{f\circ\var})X-\frac{f'\circ\var}{f}\hat{\Phi}'(\tau)\}\int_{S_o}\<H,D\rho\>d\rho d\theta\nonumber\\
&&=\int_{\S}[\Phi_1(1)-p(1)\tau^{n-1}(1)]\<H,D\rho\>d\S\nonumber\\
&&\quad-\lim_{\varepsilon\rw0}f^{n-1}(\varepsilon)(\Phi_1-p\tau^{n-1})\int_{S_o}\<H,D\rho\>d\theta,\label{46n}\eeq where $\S$ is the unit geodesic sphere centered at $o$ and $p$ is given by (\ref{43*}).

Let  a radical deformation $\u$ be an equilibrium solution where $\var$ is given by (\ref{39*}) with $\var(0)=A>0.$
The hypothesis that $\pl W(d\u)/\pl F_{ij}-p(x)(\adj d\u)_{ij}\in L^1(\B)$ for $1\leq i,\,\,j\leq n$ imply that
\be f^{n-1}(\Phi_1-p\tau^{n-1}),\quad f^{n-1}(\Phi_2-p\tau^{n-2}\var')\in L^1(0,1),\label{46*}\ee and hence
\be f^{n-1}[\Phi_2-\tau^{-n}\Phi_1]\in L^1(0,1).\label{44}\ee

Since $f(0)=0$ and $f'(0)=1,$ we take $1>\rho_0>0$ such that
$$\frac{f^n(A)}{f^n(\rho)}>\sup_{0\leq\rho\leq\rho_0}\frac{f'\circ\var(\rho)}{f'(\rho)}\qfq 0\leq\rho\leq\rho_0.$$
We have
\beq &&(n-1)\int_0^1f^{n-1}[\Phi_2-\tau^{-n}\Phi_1]d\rho=\int_0^1f^{n-1}\hat{\Phi}'(\tau)d\rho\nonumber\\
&&=\int_{\tau(\rho_0)}^\infty\frac{f^n\tau^{n-1}}{f'\tau^n-f'\circ\var}\hat{\Phi}'(\tau)d\tau+\int_{\rho_0}^1f^{n-1}\hat{\Phi}'(\tau)d\rho.\label{45}\eeq Since
$f^n=f^n\circ\var\tau^{-n},$ it follows from (\ref{45}) that the relation (\ref{44}) holds  if and only if the relation
(\ref{42}) is true.

Conversely, let  assumption (\ref{42}) hold and $p$ be given by (\ref{43*}). We  prove that
$\u$ is an equilibrium solution. Similar arguments as in the proof of Theorem 4.3 in \cite{Ball1} show that the relations (\ref{46*}) hold, that is,
$$\pl W(d\u)/\pl F_{ij}-p(x)(\adj d\u)_{ij}\in L^1(\B)\qfq 1\leq i,\,\,j\leq n.$$
Next, using (\ref{46*}) and (\ref{44n}), we deduce
$$\lim_{\varepsilon\rw0}f^n(\varepsilon)(\Phi_1-p\tau^{n-1})=0.$$  Thus, by (\ref{46n}) and (\ref{nnn2.42}), we obtain
$$\frac{\pl I(\v)}{\pl t}\Big|_{t=0}=0.$$\hfill$\Box$

\setcounter{equation}{0}
\section{Cavitation in the Incompressible Case}
\def\theequation{3.\arabic{equation}}
\hskip\parindent
Let the radial curvature $\kappa$ satisfy
\be \mu_+(\delta_0)\leq1\quad\mbox{for some}\quad\delta_0>1.\label{delta0}\ee Then estimates (\ref{kappa4}) and (\ref{kappa5*}) hold.

We need the following.

\begin{lem} \label{cnl3.1} $(i)$\,\,\,Let $1\geq\rho_1>0$ be given such that
\be\kappa(\rho)f^2(\rho)+nf'^2(\rho)>0\qfq 0\leq\rho\leq\rho_1.\label{cn3.8}\ee We fix $0<\rho_0\leq\rho_1$ such that $\si(\rho_0)<\rho_1.$
Then, for all $0< A\leq b,$ we have
$$\tau'(\rho)<0\quad\mbox{for all}\quad 0<\rho\leq\rho_0,$$ where \be b=\min\{\,\si^{-1}(\rho_1-\si(\rho_0)),\si^{-1}(\si(\delta_0)-\si(1))\}.\label{b3.2}\ee

$(ii)$\,\,\, For $0< A\leq b$ and $\tau\in[\tau(\rho_0),\infty)$ given, we solve
$f\circ\var=\tau f$ and $\si(\var)=\si(\rho)+\si(A)$ together to have $\var=\var(A,\tau)$ and $\rho=\rho(A,\tau).$ Then
\be\var_A=\frac{\tau f'f^{n-1}(A)}{f^{n-1}(f'\tau^n-f'\circ\var)},\quad\rho_A=\frac{f'\circ\var f^{n-1}(A)}{f^{n-1}(f'\tau^n-f'\circ\var)},\label{nnn3.10}\ee for
$0< A\leq b$ and $\tau\in[\tau(\rho_0),\infty).$ Moreover, there are $c_1\geq c_0>0$ such that
\be \frac{c_0A}{(\tau^n-1)^{1/n}}\leq\rho\leq\frac{c_1A}{(\tau^n-1)^{1/n}},\label{rho3.5}\ee
\be \frac{c_0A\tau}{(\tau^n-1)^{1/n}}\leq\var\leq\frac{c_1A\tau}{(\tau^n-1)^{1/n}},\label{rhov3.6}\ee
\be c_0\frac{A(\tau-1)}{(\tau^n-1)^{1/n}}\leq\var-\rho\leq c_1\frac{A(\tau-1)}{(\tau^n-1)^{1/n}},\label{rho3.6}\ee for all $(A,\tau)\in(0,b]\times[\tau(\rho_0),\infty).$

$(iii)$\,\,\,For $\tau>1$ given,
\be\lim_{A\rw0+}\frac{A}{\rho}=(\tau^n-1)^{1/n},\quad\lim_{A\rw0+}\frac{A}{\var}=\tau^{-1}(\tau^n-1)^{1/n} .\label{rhov3.8}\ee
\end{lem}

{\bf Proof}\,\,\,(i)\,\,\,Since
$$\tau'=\frac{f'\circ\var-f'\tau^n}{f\tau^{n-1}}\qfq0<\rho\leq1,$$ $\tau'<0$ if and only if
$$\frac{f'}{f^n}>\frac{f'\circ\var}{f^n\circ\var}.$$ On the other hand, condition (\ref{cn3.8}) implies
$$\Big(\frac{f'}{f^n}\Big)'=-\frac{\kappa f^2+nf'^2}{f^{n-1}}<0\qfq0<\rho\leq\rho_1.$$ Thus, $f'/f^n$ is strictly decreasing for $\rho\in(0,\rho_1].$
For $0<\rho\leq\rho_0$ and $0< A\leq\si^{-1}(\rho_1-\si(\rho_0)),$ we have
$$\rho<\var(\rho)\leq\rho_1,$$ which yields $f'/f^n>f'\circ\var/f^n\circ\var.$

(ii)\,\,\,We differentiate $f\circ\var=\tau f$ and $\si(\var)=\si(\rho)+\si(A),$ respectively, with respect to the variable $A,$ and obtain (\ref{nnn3.10}).

It follows from (\ref{39*}) that
$$\int_\rho^\var f^{n-1}ds=\int_0^Af^{n-1}ds.$$ Estimates (\ref{kappa4}) and (\ref{kappa5*}) yield
\be c_0(\var^n-\rho^n)\leq A^n\leq c_1(\var^n-\rho^n),\label{rho3.7}\ee for some $c_1\geq c_0>0.$
Since
$$\tau^n-1=\frac{1}{f^n}[f^n\circ\var-f^n]=\frac{n}{f^n}\int_\rho^\var f^{n-1}f'ds,$$ using (\ref{kappa4}) and (\ref{kappa5*}), we have
\be\frac{c_0}{\rho^n}(\var^n-\rho^n)\leq\tau^n-1\leq\frac{c_1}{\rho^n}(\var^n-\rho^n),\label{rho3.8}\ee for some $c_1\geq c_0>0.$ Thus, (\ref{rho3.5})
follow from (\ref{rho3.7}) and (\ref{rho3.8}). Similar arguments yield (\ref{rhov3.6}) and (\ref{rho3.6}).

(iii)\,\,\,By (\ref{rho3.5}), for $\tau>1$ given, $$\frac{1}{c_1}(\tau^n-1)^{1/n}\leq A/\rho\leq \frac{1}{c_0}(\tau^n-1)^{1/n}$$ for $A>0$ small. Let $a=\lim_{A\rw0+}A/\rho.$ Using (\ref{rho3.5}) and (\ref{nnn3.10}), we have
$$a=\lim_{A\rw0+}\frac{f(A)}{f(\rho)}=\lim_{A\rw0+}\frac{f'(A)}{f'(\rho)\rho_A}=\frac{\tau^n-1}{a^{n-1}},$$ which yields $a=(\tau^n-1)^{1/n}.$
A similar computation gives the second formula.
\hfill$\Box$

\begin{lem}\label{cnl3.2} Let $\rho_0\in(0,1]$ be given in Lemma $\ref{cnl3.1}$ and  let $\delta>1$ be given. Suppose $b>0$ is given such that $1<\tau(\rho_0)\leq\delta$ for all $0<A\leq b.$
Let
\be p(A,\tau)=\frac{f'\circ\var(\tau-1)}{f'\tau^n-f'\circ\var}\qfq (A,\tau)\in(0,b]\times[\tau(\rho_0),\delta].\label{cnp3.4}\ee Then
\be\min\{a_0,\,\frac{1}{n\delta}\}\leq p(A,\tau)\leq\max\{a_1,\,\frac{1}{n}\}\qfq (A,\tau)\in(0,b]\times[\tau(\rho_0),\delta],\label{cn3.4}\ee where
$$a_0=\inf_{0\leq\rho\leq\rho_0}\frac{f'^2}{\kappa f^2+nf'^2},\quad a_1=\sup_{0\leq\rho\leq\rho_0}\frac{f'^2}{\kappa f^2+nf'^2}.$$
\end{lem}

{\bf Proof}\,\,\,Let
$$a=\inf_{0<A\leq b,\,\tau(\rho_0)\leq\tau\leq\delta}p(A,\tau).$$ Suppose $A_k\rw0+$ and $\tau_k\in[\tau(\rho_0),\delta]$ such that
$$p(A_k,\tau_k)\rw a\quad\mbox{as}\quad k\rw\infty.$$ We may suppose $\tau_k\rw\tau_0$ where $\tau_0\in[1,\delta].$ We may also assume $\rho(A_k,\tau_k)\rw \hat{\rho_0}$ where $\hat{\rho}_0\in[0,\rho_0].$

Since $\tau-1=\int_\rho^\var f'(s)ds/f$ and $f'/f'\circ\var-1=\int_\rho^\var\kappa fds/f'\circ\var,$ we have
\be p(A,\tau)=\frac{1}{\frac{f\int_\rho^\var\kappa fds}{f'\circ\var\int_\rho^\var f'ds}\tau^n+\tau^{n-1}+\cdots+1}.\label{cn3.12}\ee Thus, we obtain
$$a=\frac{f'^2(\hat{\rho}_0)}{\kappa(\hat{\rho}_0)f^2(\hat{\rho}_0)\tau_0^n+(\tau_0^{n-1}+\cdots+1)f'^2(\hat{\rho}_0)}.$$
If $\tau_0>1,$ by Lemma \ref{cnl3.1}, $\hat{\rho}_0=0$ and $a\geq 1/(n\delta).$ If $\tau_0=1,$ then $a\geq a_0.$ Similar arguments yield the right hand side of (\ref{cn3.4}). \hfill$\Box$\\

Let the function $\Phi(v_1,\cdots,v_n)$ be given by (\ref{2.15}). It is said  that {\it the Baker-Ericksen
inequalities} hold if $$\frac{v_i\Phi_i-v_j\Phi_j}{v_i-v_j}\geq0\qfq i\not=j,\quad v_i\not=v_j.$$
We further assume that $W$ is bounded below. By Theorem \ref{t2},  the radical deformation $\var,$ given by (\ref{39*}), with $\var(0)=A>0$, is an equilibrium solution
if and only if
\be\frac{v^{n-1}}{(v^n-1)^2}\hat{\Phi}'(v)\in L^1(\delta,\infty)\qfq \delta>1.\label{ab3.4}\ee
 Let $p$ be given by (\ref{43*}). Then the radical component of the Cauchy stress tensor is given by
$$T=\tau^{1-n}\Phi_1-p\qfq\rho>0.$$
Let $\rho_0\in(0,1]$ be given in Lemma \ref{cnl3.1}. It follows from (\ref{43*}) and Lemma \ref{cnl3.1} (ii) that $T(0)=\lim_{\rho\rw0+}T$ exists if and only if
the integral
\be \int_{\tau(\rho_0)}^\infty\frac{f'\circ\var}{f'\tau^n-f'\circ\var}\hat{\Phi}'(\tau)d\tau\ee converges.

The total stored energy of the deformation is given by
\be E(A)=\omega_n\int_0^1f^{n-1}(\rho)\Phi(\rho)d\rho,\label{x3.6}\ee where $\omega_n$ is the area of the unit sphere $S_o$ in $M_o$ and $\Phi(\rho)=\Phi(\var'(\rho),\tau(\rho),\cdots,\tau(\rho)).$
We define $E(0)=\omega_n\si(1)\Phi(1).$

\begin{pro}  Let $(M,g,o)$ be a model with $\mu_+(1)\leq1$ and let $(\ref{ab3.4})$ hold. Let $(\ref{39*})$ be an equilibrium solution with $A>0$ and
$\mu_+(\var(1))\leq1.$ Then

$(i)$\,\,\,Then $T(0)$ exists and is finite if and only if $E(A)<\infty.$

$(ii)$\,\,\,  Let $\Phi$ satisfy the Baker-Ericksen inequalities.Then $T$ is
an increasing function in $\rho>0.$
\end{pro}

{\bf Proof}\,\,\,Using equation (\ref{44n}), we have
$$T'(\rho)=\frac{f'\circ\var}{f}\tau^{1-n}\hat{\Phi}'(\tau)=(n-1)\frac{f'\circ\var}{f\tau^{2n-1}}(\tau^n-1)\frac{\tau\Phi_2-\tau^{1-n}\Phi_1}{\tau-\tau^{1-n}}\geq0,$$
that is, (ii) is true.
Using the formula above, we obtain
$$[f^n\Phi]'=nf'f^{n-1}\Phi+f^n\hat{\Phi}'(\tau)\tau'=nf'f^{n-1}\Phi+(f^n-\frac{f'f^n\circ\var}{f'\circ\var})T'.$$ Thus,
\beq n\int_\rho^1f'f^{n-1}\Phi d\rho&&=f^n(1)\Phi(1)-f^n(\rho)\Phi(\rho)+\int_\rho^1(\frac{f'f^n\circ\var}{f'\circ\var}-f^n)T'd\rho\nonumber\\
&&=f^n(1)\Phi(1)-f^n(\rho)\Phi(\rho)+(\frac{f'f^n\circ\var}{f'\circ\var}-f^n)T\Big|_{\rho=1}\nonumber\\
&&\quad-(\frac{f'f^n\circ\var}{f'\circ\var}-f^n)T+\int_{\rho}^1(\frac{f'f^n\circ\var}{f'\circ\var}-f^n)_\rho Td\rho,\label{n3.1}\eeq for $0\leq\rho\leq1.$
Using (\ref{nn2.12}), we have
\be\lim_{\rho\rw0+}\frac{1}{f(\rho)}(\frac{f'f^n\circ\var}{f'\circ\var}-f^n)_\rho=-\kappa(o)f^n(A).\label{cnn3.8}\ee In addition, it follows from (\ref{43*}) that
\beq|fT|&&\leq\int_\rho^1f'\circ\var\tau^{1-n}|\hat{\Phi}'(\tau)|d\rho=\int_{\rho_0}^1f'\circ\var\tau^{1-n}|\hat{\Phi}'(\tau)|d\rho\nonumber\\
&&\quad+\int_{\tau(\rho_0)}^\infty \frac{f\circ\var f'\circ\var}{\tau (f'\tau^{n}-f'\circ\var)}|\hat{\Phi}'(\tau)|d\tau.\nonumber\eeq
Thus, by (\ref{ab3.4}),  $f|T|$ is bounded above. Therefore the last integral in the right hand side of (\ref{n3.1}) converges by (\ref{cnn3.8}). Moreover, $-f^n\Phi$ is also bounded above
since $\Phi$ is bounded below. Thus, the left hand side of (\ref{n3.1}) exists and is finite if and only if $T(0)$ exits and is finite. Then (i) follows by (\ref{kappa5*}).
 \hfill$\Box$\\

We consider the case when there is a force acting on the unit geodesic sphere $\S.$ A function on $M$ is said to be a {\it force density}. We say that $\u\in W^{1,1}(\B,M)$ is an {\it equilibrium solution} to the boundary value problem with corresponding pressure $p$ if $\det d\u=1$ a.e. in $\B,$
$$\pl W(d\u)/\pl F_{ij}-p(x)(\adj d\u)_{ij}\in L^1(\B)\qfq 1\leq i,\,\,j\leq n,$$ and
$$\int_\B\<D_FW-pD_F\det d\u,\,\,DH\>dg-\int_\S\<Dq,H\>d\S=0\qfq H\in C^\infty(\bar{\B},\Ga(\zeta)),$$ where $q$ is a given force density which is a differentiable function on $M.$ They are the Euler-Lagrange equations for the functional
$$I_1(\u)=\int_\B\{W(d\u)-p(\det d\u-1)\}dg-\int_\S q(\u)d\S.$$
A force density $q$ is said to be {\it radical with respect to $o$} if there is $\hat{q}\in C^1(\R)$ such that
$$q(x)=\hat{q}(\rho(x))\qfq x\in M.$$ In particular, we take $$q(x)=P\rho(x)\qfq x\in M,$$ where $P\in\R$ is a constant.

It is easy to check that the identity deformation $\var(\rho)=\rho$ is an equilibrium solution to the boundary value problem
with the corresponding pressure
$$p=\Phi_1(1)-P.$$

Let all the assumptions in Theorem \ref{t2} hold. Let $\var$ be given by (\ref{39*}) with $A>0.$ By similar arguments  as in the proof of Theorem \ref{t2} that $\var$ is an equilibrium
 solution to the above boundary value problem if and only if it is an equilibrium solution with the corresponding
pressure
$$p=\tau^{1-n}\Phi_1-T,$$ where
\be T=\frac{P}{\tau^{n-1}(1)}-\int^{1}_{\rho}
\frac{f'\circ\var}{f\circ\var}\tau^{2-n}\hat{\Phi}'(\tau)d\rho\label{2.1}\ee for $\rho>0$ and $\lam>0.$ As a natural boundary condition given in \cite{Ball1} to get a unique solution for the corresponding $\lam\in\R,$ we assume that
\be T(0)=0.\label{2.2}\ee
The total energy of the deformation (\ref{39*}) with $\var(0)=A\geq0$ is given by
\beq I(A)=\int_\B W(d\u)dg-P\int_{\S}\rho(\u)d\S=E(A)-\omega_nf^{n-1}(1)P\var(1),\label{xn3.21}\eeq where $E(A)$ is given by (\ref{x3.6}).

By (\ref{2.1})  the possible values of $A$ such that (\ref{2.2}) is satisfied are the roots of the equation
\be P=\chi(A)\label{rh3.22}\ee where \be \chi(A)=\tau^{n-1}(1)\int_0^1\frac{f'\circ\var}{f\circ\var}\tau^{2-n}\hat{\Phi}'(\tau)d\rho.\label{2.5}\ee
Thus, the bifurcation from the trivial solution is governed by the behavior of $\chi(A)$ as $A\rw0+.$

Let $$\hat{I}(A)=\int_0^1f^{n-1}(\rho)\Phi(\rho)d\rho.$$ Let $\rho_0\in(0,1]$ be given by Lemma \ref{cnl3.1}. Since
\beq\hat{I}(A)&&=\int_{\rho_0}^1f^{n-1}(\rho)\Phi(\rho)d\rho+\si(\rho_0)\hat{\Phi}(\tau(\rho_0))-\int_0^{\rho_0}\si(\rho)\hat{\Phi}'(\tau)\tau'd\rho\nonumber\\
&&=\int_{\rho_0}^1f^{n-1}(\rho)\Phi(\rho)d\rho+\si(\rho_0)\hat{\Phi}(\tau(\rho_0))+\int_{\tau(\rho_0)}^{\infty}\si(\rho(A,\tau))\hat{\Phi}'(\tau)d\tau,\label{xn3.24}\eeq using the formula
in (\ref{nnn3.10})  we have
\beq \hat{I}'(A)&&=\int_{\rho_0}^1f^{n-1}(\rho)\hat{\Phi}'(\tau(\rho))\tau_A(\rho)d\rho+\int_{\tau(\rho_0)}^\infty f^{n-1}\rho_A\hat{\Phi}'(\tau)d\tau\nonumber\\
&&=f^{n-1}(A)\Big[\int_{\rho_0}^1\frac{f'\circ\var}{f\circ\var}\tau^{2-n}\hat{\Phi}'(\tau)d\rho
+\int_{\tau(\rho_0)}^\infty\frac{f'\circ\var}{f'\tau^n-f'\circ\var}\hat{\Phi}'(\tau)d\tau\Big]\nonumber\\
&&=f^{n-1}(A)\int_{0}^1\frac{f'\circ\var}{f\circ\var}\tau^{2-n}\hat{\Phi}'(\tau)d\rho\nonumber\eeq which yields
\be I'(A)=\omega_nf^{n-1}(A)\tau^{1-n}(1)[\chi(A)-P]\qfq A\geq0.\label{rh3.24}\ee

We suppose
\be\frac{\hat{\Phi}'(\tau)}{\tau^n-1}\in L^1(\delta,\infty)\qfq\delta>1.\label{cn3.5}\ee
Let $\hat{\Phi}(v)$ be twice differentiable at $v=1.$ Thus, (\ref{cn3.5}) and $\hat{\Phi}'(1)=0$ imply
$$\frac{\hat{\Phi}'(\tau)}{\tau^n-1}\in L^1(1,\infty).$$ Using (\ref{rho3.5}), we have
$$\Big|\int_{\tau(\rho_0)}^{\infty}\si(\rho(A,\tau))\hat{\Phi}'(\tau)d\tau\Big|\leq cA^n\int_1^\infty\frac{1}{\tau^n-1}|\hat{\Phi}'(\tau)|d\tau.$$
By (\ref{xn3.24}) and (\ref{xn3.21}), we obtain
\be \lim_{A\rw0+}I(A)=\omega_n[\si(1)\Phi(1)-f^{n-1}(1)P]=I(0).\label{xn3.27}\ee

\begin{lem}\label{cnl3.3} Let $\chi$ be given by $(\ref{2.5})$ and let $\hat{\Phi}(v)$ be twice differentiable at $v=1.$ Then
\be\lim_{A\rw0+}\chi(A)=\int_1^\infty\frac{1}{\tau^n-1}\hat{\Phi}'(\tau)d\tau.\label{cn3.10}\ee
\end{lem}

{\bf Proof}\,\,\,Let $\rho_0\in(0,1]$ be given in Lemma \ref{cnl3.1}. We have
\be\int_0^1\frac{f'\circ\var}{f\circ\var}\tau^{2-n}\hat{\Phi}'(\tau)d\rho=\int_{\rho_0}^1\frac{f'\circ\var}{f\circ\var}\tau^{2-n}\hat{\Phi}'(\tau)d\rho
+\int_{\tau(\rho_0)}^\infty\frac{f'\circ\var}{f'\tau^n-f'\circ\var}\hat{\Phi}'(\tau)d\tau,\label{cn3.11}\ee for  $0< A\leq\si^{-1}(\rho_1-\si(\rho_0)).$
Since $\Phi'(1)=0,$ the first integral in the right hand side of (\ref{cn3.11}) goes to zero as $A\rw0+.$

Since
$$\hat{\Phi}'(\tau)/(\tau-1)\rw \Phi''(1)\quad \mbox{as}\quad \tau\rw0+,$$ we fixed $\delta_0\geq\delta>1$ such that $|\hat{\Phi}'(\tau)|/(\tau-1)\leq |\hat{\Phi}''(1)|+1$
for all $1\leq\tau\leq\delta.$ By Lemma \ref{cnl3.2},
$$\int_{\tau(\rho_0)}^\delta\frac{f'\circ\var}{f'\tau^n-f'\circ\var}|\hat{\Phi}'(\tau)|d\tau\leq c[|\hat{\Phi}''(1)+1],$$ for $A>0$ small. Thus, (\ref{cn3.10}) follows from the dominated
convergence theorem. \hfill$\Box$\\

Let
$$P_{cr}=\int_1^\infty\frac{1}{\tau^n-1}\hat{\Phi}'(\tau)d\tau.$$ The physical meaning of $P_{cr}$ is given in \cite{GeLi}, \cite{Ball1}: (\ref{rh3.24}), (\ref{xn3.27}) and Lemma \ref{cnl3.3} show that the trivial solution $A=0$ is a local minimum (resp. local maximum) if $P<P_{cr}$ (resp. $P>P_{cr}$).  The following proposition illustrates
the close relations among the radial curvatures, the constitutive function $\hat{\Phi}(v),$ and the behavior  of $\chi(A)$ as $A\rw0+.$

\begin{pro} \label{inp3.2} Let $\hat{\Phi}(v)$ be twice differentiable at $v=1.$ Then $\chi'(0)=0.$ Furthermore, the following holds.

$(i)$\,\,\, Let $n=2.$ Let $\Phi$ satisfy the Baker-Ericksen inequalities with $\hat{\Phi}''(1)>0.$ Suppose there is $\varepsilon\in(0,1]$ such that the radial curvature $\kappa$ is a constant $\kappa_0$ for all $\rho\in(0,\varepsilon].$ If $\kappa_0=0,$ then
\be P_{cr}-\frac{1}{2}[1+\frac{f^2(1)}{f'(1)}\int_\varepsilon^1\frac{\kappa}{f}d\rho]\hat{\Phi}''(1)>0\quad \mbox{$($ resp. $<0$ $)$}\ee implies $\chi'(A)>0$ $($ resp. $<0)$ for $A>0$ small. If $\kappa_0\not=0,$ then
\be \kappa_0>0\quad \mbox{$($ resp. $<0$ $)$}\ee implies $\chi'(A)<0$ $($ resp. $>0)$ for $A>0$ small.

$(ii)$\,\,\,Let $n\geq3.$ Let $\kappa(o)\not=0.$ Then
\be\kappa(o)\int_1^\infty\frac{(\tau^2-1)\tau^n}{(\tau^n-1)^{2(1+1/n)}}\hat{\Phi}'(\tau)d\tau<0\quad \mbox{$($ resp. $>0$ $)$}\label{rho3.30*}\ee
implies $\chi'(A)>0$ $($ resp. $<0)$ for $A>0$ small. In addition,  if there is some $\varepsilon\in(0,1]$ such that $\kappa=0$ for $\rho\in(0,\varepsilon),$ then
\be P_{cr}-\frac{1}{n(n-1)}\Big[1+\frac{f^n(1)}{f'(1)}\int_\varepsilon^1\frac{\kappa}{f^{n-1}}d\rho\Big]\hat{\Phi}''(1)>0\quad \mbox{$($ resp. $<0$ $)$}\label{rhov3.31}\ee
implies $\chi'(A)>0$ $($ resp. $<0)$ for $A>0$ small.
\end{pro}

{\bf Proof}\,\,\,Let $\rho_0\in(0,1]$ be given in Lemma \ref{cnl3.1}. Using (\ref{cn3.11}) and (\ref{nnn3.10}), we have
\beq \chi'(A)&&=(n-1)\frac{\tau_A(1)}{\tau(1)}\chi(A)+\tau^{n-1}(1)[\int_{\rho_0}^1(\frac{f'\circ\var}{f\circ\var}\tau^{2-n})_A\hat{\Phi}'(\tau)d\rho\nonumber\\
&&\quad
+\int_{\rho_0}^1I_1\hat{\Phi}''(\tau)\tau_Ad\rho]+
\tau^{n-1}(1)[-p(A,\tau)\frac{\hat{\Phi}'(\tau)}{\tau-1}\Big|_{\tau=\tau(\rho_0)}\tau_A(\rho_0)\nonumber\\
&&\quad+f^{n-1}(A)\int_{\tau(\rho_0)}^\infty
I_2\hat{\Phi}'(\tau)d\tau],\label{cn3.23}\eeq
where
$$I_1=\frac{f'\circ\var}{f\circ\var}\tau^{2-n},\quad I_2=\frac{(\kappa f'^2\circ\var-\kappa\circ\var f'^2\tau^2)\tau^n}{f^{n-2}(f'\tau^n-f'\circ\var)^3},$$ and
$p(A,\tau)$ is given by (\ref{cnp3.4}). Clearly, all the terms in the right hand side of (\ref{cn3.23}) go to zero if the last term converges to zero as $A\rw0+.$

Let $\delta>1$ be given. It follows from (\ref{rho3.5}) that
\be A\leq\frac{\rho_0}{c_0}(\tau^n-1)^{1/n}\leq c(\tau-1)^{1/n}\qfq (A,\tau)\in(0,b)\times[\tau(\rho_0),\delta].\label{rho3.30}\ee

Using (\ref{rho3.5}), (\ref{rho3.6}), (\ref{cn3.4}), and (\ref{rho3.30}), we have
\beq |I_2|&&=\Big|\frac{-f'^2\circ\var\int_\rho^\var\kappa'(s)ds+\kappa\circ\var[f'^2\circ\var(1-\tau^2)+2\tau^2\int_\rho^\var\kappa ff'ds]}{f^{n-2}f'^3\circ\var(\tau-1)^3}\Big|p^3\tau^n\nonumber\\
&&\leq
c\frac{A(\tau-1)^{1-1/n}+\tau-1}{A^{n-2}(\tau-1)^{2(1+1/n)}}\leq \frac{c}{A^{n-2}(\tau-1)^{1+2/n}}\qfq (A,\tau)\in(0,b)\times[\tau(\rho_0),\delta].\nonumber\eeq
Thus, by (\ref{rho3.30}) again,
\be f^{n-1}(A)\Big|I_2\hat{\Phi}'(\tau)\Big|\leq c\frac{A^{1-\a}}{(\tau-1)^{(2-\a)/n}}\Big|\frac{\hat{\Phi}'(\tau)}{\tau-1}\Big|\qfq (A,\tau)\in(0,b)\times[\tau(\rho_0),\delta] \label{rho3.31}\ee
and for $0\in[0,1].$ It follows from (\ref{cn3.23}) and (\ref{rho3.31}) that $\chi'(0)=\lim_{A\rw0+}\chi'(A)=0.$

(i)\,\,\,We may assume that $0<\rho_0<\varepsilon$ is small enough such that $\var\leq\varepsilon$ when $\rho\in(0,\rho_0].$ Let $\kappa_0=0.$ Then $f'(\rho_0)=1$, $f(\rho_0)=\rho_0,$ and $I_2=0.$
In this case we have, by (\ref{cn3.23}) and (\ref{cn3.12}),
\beq \lim_{A\rw0+}\frac{\chi'(A)}{f(A)}&&=\frac{f'(1)}{f^2(1)}P_{cr}+\lim_{A\rw0+}\int_{\rho_0}^1\frac{f'\circ\var}{ff\circ\var}I_1\hat{\Phi}''(\tau)d\rho-\lim_{A\rw0+}\frac{p(A,\tau)f'\circ\var}{ff\circ\var}
\frac{\hat{\Phi}'(\tau)}{\tau-1}\Big|_{\tau=\tau(\rho_0)}\nonumber\\
&&=\frac{f'(1)}{f^2(1)}P_{cr}+\Big(\int_{\rho_0}^1\frac{f'^2}{f^3}d\rho-\frac{1}{2\rho_0^2}\Big)\hat{\Phi}''(1)\nonumber\\
&&=\frac{f'(1)}{f^2(1)}P_{cr}-\frac{1}{2}\Big[\frac{f'(1)}{f^2(1)}+\int^1_{\varepsilon}\frac{\kappa}{f}d\rho\Big]\hat{\Phi}''(1).\nonumber\eeq
Thus, the case $\kappa_0=0$ follows.

Let $\kappa_0\not=0.$ Let $q=f'/f\circ\var.$ By $(\ref{cn3.8}),$ $(f'/f)'=-(\kappa f^2+f'^2)/f^2<0$ for $\rho\in(0,\rho_0].$ Thus,
$$q^2\tau^2>1\qfq (A,\tau)\in(0,b]\times[\tau(\rho_0),\infty).$$
Let $\delta>1$ be fixed such that
$$\frac{1}{\tau-1}\hat{\Phi}'(\tau)\geq\frac{1}{2}\hat{\Phi}''(1)\qfq\tau\in[1,\delta].$$ We have
$$\Big|\frac{(q^2-1)\tau^2}{\tau-1}+\tau+1\Big|\geq2-\Big|\frac{2\kappa_0f\int_\rho^\var ff'ds}{f'^2\circ\var\int_\rho^\var f'ds}\Big|\tau^2
\geq 2-c|\kappa_0|f^2(\rho_0),$$ for $(A,\tau)\in(0,b]\times[\tau(\rho_0),\delta].$ We assume that $\rho_0\in(0,\varepsilon]$ is also such that
$2-c|\kappa_0|f^2(\rho_0)>0.$ Thus, by (\ref{cn3.4}), we obtain
\beq -\frac{1}{\kappa_0}I_2\hat{\Phi}'&&=\frac{(q^2\tau^2-1)\tau^2p(A,\tau)}{f'\circ\var(\tau-1)^3}\hat{\Phi}'(\tau)\geq\frac{c_0}{\tau-1}\hat{\Phi}''(1)\qfq(A,\tau)\in(0,b]\times[\tau(\rho_0),\delta],\nonumber\eeq  which yields, by (\ref{cn3.23}),
$$\lim_{A\rw0+}\frac{\chi'(A)}{f(A)}=\left\{\begin{array}{l}-\infty\quad\mbox{if}\quad\kappa_0>0;\\
+\infty\quad\mbox{if}\quad\kappa_0<0.\end{array}\right.$$ Thus, the case $\kappa_0\not=0$ follows.

(ii)\,\,\,Since $n\geq3,$ (\ref{rho3.31}) with $\a=0$ and (\ref{cn3.5}) imply the integral
$$\int_{\tau(\rho_0)}f^{n-2}(A)I_2\hat{\Phi}'(\tau)d\tau$$ converges as $A\rw0+.$ Using (\ref{cn3.23}), (\ref{rho3.5}), (\ref{rhov3.6}), and (\ref{rhov3.8}), we obtain
$$\lim_{A\rw0+}\frac{\chi'(A)}{f(A)}=\kappa(o)\int_1^\infty\frac{(1-\tau^2)\tau^n}{(\tau^n-1)^{2(1+1/n)}}d\tau,$$ which gives (\ref{rho3.30*}). Finally, a similar computation as in (i) for the case
$\kappa_0=0$ yields (\ref{rhov3.31}). \hfill$\Box$

\begin{rem}
If $\Phi$ satisfies the Baker-Ericksen inequalities, then $(\ref{rho3.30*})$ is equivalent to
$$\kappa(o)<0\quad \mbox{$($ resp. $>0$ $)$}.$$
\end{rem}

Let $P>0.$ If $A_0$ is a root of (\ref{rh3.22}), that is,
$$P=\chi(A_0),$$ then from (\ref{rh3.24})
$$I''(A_0)=\omega_nf^{n-1}(A_0)\tau^{1-n}(1)\chi'(A_0).$$ Let $\Phi$ satisfy the Baker-Ericksen inequalities and let $n\geq3.$ From (\ref{rho3.30*}), $A_0$ is a local minimum (resp. local maximum)
of $I$ if $\kappa(o)<0$ (resp. $\kappa(o)>0$) (for $A_0$ small).

\setcounter{equation}{0}
\section{Cavitation in the Compressible Case}
\def\theequation{4.\arabic{equation}}
\hskip\parindent By Theorem \ref{nt1.1}, an equilibrium solution $\var$ satisfies the equation
\be [f^{n-1}(\rho)\Phi_1(\rho)]_\rho=(n-1)f^{n-2}(\rho)f'\circ\var(\rho)\Phi_2(\rho)\qfq x\in\Om,\quad
\rho(x)>0,\label{3.1}\ee where
$$\Phi_1(\rho)=\Phi_1(\var',\tau,\cdots,\tau),\quad \Phi_2(\rho)=\Phi_2(\var',\tau,\cdots,\tau),
\quad\tau(\rho)=\frac{f\circ\var(\rho)}{f(\rho)}.$$

Let $\var$ be a solution to problem (\ref{3.1}) with $\var'>0$ and $\var>0$ on $(0,1].$ We define
\be T(\rho)=\tau^{n-1}(\rho)\Phi_1(\rho)\qfq\rho\in(0,1],\label{3.2*}\ee which is the radial component of the Cauchy stress (\cite{Ball1}). By (\ref{3.1}),
\be T'(\rho)=(n-1)\frac{f'\circ\var}{f\circ\var}\tau^{1-n}(\tau\Phi_2-\var'\Phi_1)\qfq\rho\in(0,1].\label{3.3*}\ee

It follows from (\ref{3.3*}) that

\begin{pro}\label{p3.1}Let $\var$ be a solution of $(\ref{3.1})$ with $\var'(\rho)>0$ for all $\rho\in(0,1].$ If the Baker-Ericksen inequalities hold,
then
\be T'(\rho)[\var'(\rho)-\tau(\rho)]\leq 0\qfq \rho\in(0,1].\label{3.4*}\ee
\end{pro}

Let
\be \tilde{T}(\rho)=\Phi(\rho)-\var'(\rho)\Phi_1(\rho),\label{3.5*}\ee where $\Phi(\rho)=\Phi(\var',\tau,\cdots,\tau).$ $\tilde{T}$ is said to be
the radial component of the inverse Cauchy stress (\cite{Ball1}). We obtain by (\ref{3.1}) and (\ref{3.3*})
\be\tilde{T}'(\rho)=-\frac{f'}{f'\circ\var}\tau^nT'\qfq\rho\in(0,1].\label{3.6*}\ee
It follows (\ref{3.6*}), (\ref{3.3*}) and (\ref{3.1}) that
\beq &&\Big\{f^n[\Phi-(\var'-\tau)\Phi_1]\Big\}'=[f^n\tilde{T}+f\circ\var f^{n-1}\Phi_1]'=nf'f^{n-1}\tilde{T}+f^n\tilde{T}'\nonumber\\
&&\quad+f'\circ\var\var'f^{n-1}\Phi_1+f\circ\var(f^{n-1}\Phi_1)'\nonumber\\
&&=nf'f^{n-1}\Phi+(f'\circ\var-f')f^{n-1}[\var'\Phi_1+(n-1)\tau\Phi_2].\label{3.7*}\eeq
If $\kappa=0$ for $\rho\in(0,1],$ then $f'=1,$ $f=\rho,$ and (\ref{3.7*}) becomes
$$\Big\{\rho^n[\Phi-(\var'-\tau)\Phi_1]\Big\}'=n\rho^{n-1}\Phi,\quad \tau=\frac{\var}{\rho},$$ which is the radial version of the
conservation law (\cite{Ball1}).

\subsection{Constitutive Assumptions}
\hskip\parindent Throughout this paper unless otherwise stated we assume the  class of constitutive functions $W(F)=\Phi(v_1,\cdots,v_n)$ have the form ($n\geq2$)
\be \Phi(v_1,\cdots,v_n)=\sum_{i=1}^n\phi(v_i)+h(v_1\cdots v_n),\label{3.2}\ee where functions $\phi$ and $h$ satisfy
the following assumptions:

(A 1)\,\,\,$h:$ $(0,\infty)\rw\R$ is $C^2$ and strictly convex;

(A 2)\,\,\,$\lim_{v\rw0+}h(v)=\lim_{v\rw\infty}\dfrac{h(v)}{v}=+\infty;$

(A 3)\,\,\,$h$ satisfies
$$\underline{\lim}_{v\rw\infty}\frac{vh'(v)}{h(v)}>1,$$ and
let
$$\theta(s)=\underline{\lim}_{v\rw\infty}\frac{h(sv)}{h(v)}\qfq s\in(0,\infty),$$
and we assume that $\theta:$ $(0,\infty)\rw(0,\infty)$ is continuous;

(A 4)\,\,\,$\phi:$ $(0,\infty)\rw(0,\infty)$ is $C^2$ and convex;

(A 5)\,\,\,$v\phi'(v)$ is increasing on $(0,\infty);$

(A 6)\,\,\,Let $t_0\geq0$ be such that $\phi'(t_0)=0.$ Let
$$q_1(s)=\sup_{v>t_0}\frac{\phi'(v)}{\phi'(sv)}\qfq s>1;
\quad q_0(s)=\inf_{v>t_0/s}\frac{\phi'(v)}{\phi'(sv)}\qfq s\in(0,1].$$ We assume that $q_1\in C^1[1,\infty)$ and $q_0\in C^1(0,1]$ satisfy
\be\lim_{s\rw\infty}q_1(s)=0,\quad \lim_{s\rw0+}q_0(s)=\infty,\label{3.3}\ee
\be q_1'(s)<0\qfq s\in[1,\infty)\quad\mbox{and}\quad q_0'(s)<0\qfq s\in(0,1],\label{re3.41}\ee respectively.

(A 7)\,\,\,there are $\delta_0>0$ and  $\delta_1>0$ such that if $|s-1|<\delta_0$ then
$$|\phi'(sv)|\leq\delta_1\frac{\phi(v)}{v}\quad\mbox{for all}\quad v>0;$$

(A 8)\,\,\,$\phi(v)\leq\delta_2(1+v^\a+v^{-\b})$ for all $v>0,$ where $\delta_2>0,$ $0<\a<n,$ and $0\leq\b<1+1/(n-1);$

(A 9)\,\,\,$\lim_{v\rw\infty}\phi(v)=+\infty.$

\begin{rem} $(A \,1),$  $(A \,2),$ $(A \,4),$ $(A \,5),$ $(A \,7),$ $(A \,8),$ and $(A \,9)$ are given in $\cite{Ball1}.$ In addition, it is easy to check that $q_1$ and $q_0$ are decreasing on $[1,\infty)$ and $(0,1],$ respectively, and $q_1(1)=q_0(1)=1.$
\end{rem}

\begin{exl} Let
$$\phi(v)=\mu (v^\a-n)+\frac{\nu}{v^\b},$$ where $\mu>0,$ $\nu\geq0,$ $1<\a<n,$ and $0\leq\b<1+1/(n-1).$ Let $$h(v)=H(v)-n,$$
 where $H:$ $(0,\infty)\rw\R^+$ is a $C^3$ function and satisfies $\lim_{\delta\rw0+}H(\delta)=\infty,$ $H''(\delta)>0$ for all $\delta>0,$ and
$$H(\delta)=k(\delta-1-k^{-1})^2\qfq\delta\geq1/2,$$ where $k>1.$

Clearly, $(A\,4),$ $(A\,5),$ $(A\,7),$ $(A\,8),$ and $(A\,9)$ are true for $\phi.$ We check $(A\,6).$ Since
$$\frac{\phi'(v)}{\phi'(sv)}=\frac{h_1(v)}{h_2(v)},$$ where
$$h_1(v)=\mu\a-\frac{\b\nu}{v^{\a+\b}},\quad h_2(v)=\mu\a s^{\a-1} -\frac{\b\nu}{s^{\b+1}v^{\a+\b}},$$
$$h_1'(v)h_2(v)-h_2'(v)h_1(v)=\frac{\b(\a+\b)\nu\mu\a}{v^{\a+\b+1}}(s^{\a-1}-\frac{1}{s^{\a+\b}}),$$ we have
$$q_1(s)=\lim_{v\rw\infty}\frac{h_1(v)}{h_2(sv)}=\frac{1}{s^{\a-1}}\qfq s\in[1,\infty);\quad q_0(s)=\frac{1}{s^{\a-1}}\qfq s\in(0,1].$$ Thus, $(A\,6)$ holds.

$(A\,1)$ and $(A\,2)$ automatically hold for $h.$  $(A\,3)$ is also true since
$$\underline{\lim}_{v\rw\infty}\frac{vh'(v)}{h(v)}=2,\quad \theta(s)=s^2\qfq s\in(0,\infty).$$
\end{exl}

\subsection{Equilibrium Solutions}
\hskip\parindent
Let $\lam>0$ be given. It follows from (\ref{3.1}) and (\ref{3.2}) that  equilibrium solutions are given by
problem
\be\left\{\begin{array}{l}\var(1)=\lam,\\
f[\phi''(\var')+h''(\var'\tau^{n-1})\tau^{2(n-1)}]\var''=(n-1)[f'\circ\var\phi'(\tau)
-f'\phi'(\var')]\\
\quad-(n-1)(f'\circ\var\var'-f'\tau)h''(\var'\tau^{n-1})\var'\tau^{2n-3}\qfq \rho\in(0,1).
\end{array}\right.\label{3.4}\ee

{\bf Regular Equilibrium Solutions}\,\,\,An equilibrium solution $\var\in C^1(0,1]$ to problem (\ref{3.4}) is said to be {\it regular} if $\var(0)=\lim_{\rho\rw0+}\var(\rho)=0.$
In the case of the Euclidean space there is a unique regular equilibrium solution $\var=\lam\rho$ which plays an important role in the analysis (\cite{Ball1,Si}).
We derive some properties of regular equilibrium solutions.

We assume that $\lam>0$ is given such that the radial curvature satisfies
\be\mu_+(\max\{\lam,1\})=\int_0^{\max\{\lam,1\}}s\kappa_+(s)ds\leq1.\label{es}\ee
By Proposition \ref{p1.1}
$$f'(\rho)>0\qfq\rho\in[0,\max\{\lam,1\}].$$

Let $$b_0(\rho)=\min_{0\leq s\leq\rho}f'(s),\quad b_1(\rho)=\sup_{0\leq s\leq\rho}f'(s)\qfq\rho\in[0,\max\{\lam,1\}].$$

Set
$$\a_1(\rho,s)=\max\Big\{\,b_1(\rho)/b_0(s),\,\,q_1^{-1}\Big(b_0(\rho)/b_1(s)\Big)\,\Big\},$$
$$\a_0(\rho,s)=\min\Big\{\,b_0(\rho)/b_1(s),\,\,q_0^{-1}\Big(b_1(\rho)/b_0(s)\Big)\Big\}$$
for $(\rho,s)\in[0,1]\times[0,\lam],$ where $q_0$ and $q_1$ are given by (A 6).

Then the following lemma is immediate.

\begin{lem}\label{l3.1} Let $\mu_+(\max\{\lam,1\})\leq1$ and let $(\rho,s)\in[0,1]\times[0,\lam]$ be given. $\a_1(\cdot,s)$ and $\a_1(\rho,\cdot)$ are increasing.
$\a_0(\cdot,s)$ and $\a_0(\rho,\cdot)$ are
decreasing.
\end{lem}

We have

\begin{thm} \label{t3.2}Let $\mu_+(\max\{\lam,1\})\leq1.$ Let $(A 1),$ $(A 4),$ $(A 5),$ and $(A 6)$ hold.
If $\var\in C^1(0,1]$ is a regular equilibrium solution to problem $(\ref{3.4}),$ then
\be \a_0(\rho,\var(\rho))\tau(\rho)\leq\var'(\rho)\leq\a_1(\rho,\var(\rho))\tau(\rho)\qfq\rho\in(0,1].\label{re}\ee
\end{thm}

{\bf Proof}\,\,\,First, we prove that the right hand side of the inequalities (\ref{re}) holds true. We suppose for a contradiction that there is $\rho_0\in(0,1]$
such that
\be\var'(\rho_0)>\a_1(\rho_0,\var(\rho_0))\tau(\rho_0).\label{re1}\ee

Since $\a_1(\rho_0,\var(\rho_0))\geq b_1(\rho_0)/b_0\circ\var(\rho_0),$ by (\ref{re1}) we have
\be f\tau'= f'\circ\var\var'-f'\tau>(b_0\a_1-b_1)\tau\geq0\quad\mbox{at}\quad\rho=\rho_0.\label{re2}\ee
Next, we claim
\be f'\circ\var\phi'(\tau)-f'\phi'(\var')\leq0\quad\mbox{at}\quad\rho=\rho_0.\label{re3}\ee
We assume that $\phi'(\tau(\rho_0))\leq0.$ Since $\phi'$ is increasing, $\a_1(\rho_0,\var(\rho_0))f'\circ\var(\rho_0)\geq f'(\rho_0),$
and $\a_1(\rho_0,\var(\rho_0))\geq1,$
it follows from (\ref{re1}) and (A 5) that
\beq && f'\circ\var\phi'(\tau)-f'\phi'(\var')\leq f'\circ\var\phi'(\tau)-f'\phi'(\a_1\tau)\nonumber\\
&&=\Big[f'\circ\var-\frac{f'}{\a_1}\Big]\phi'(\tau)+\frac{f'}{\a_1\tau}[\tau\phi'(\tau)-\a_1\tau\phi'(\a_1\tau)]\leq0\quad
\mbox{at}\quad\rho=\rho_0.\nonumber\eeq Let $\phi'(\tau)>0.$ Let $$t_1=q_1^{-1}\Big(\frac{b_0(\rho_0)}{b_1\circ\var(\rho_0)}\Big).$$
Then $q_1(t_1)=b_0(\rho_0)/b_1\circ\var(\rho_0),$ that is,
$$b_1\circ\var(\rho_0)\phi'(\tau(\rho_0))\leq b_0(\rho_0)\phi'(t_1\tau(\rho_0)),$$ which implies that
(\ref{re3}) holds true since $\a_1(\rho_0,\var(\rho_0))\geq t_1.$

From (\ref{re3}), (\ref{re2}), and (\ref{3.4}) we obtain
\be\var''(\rho_0)<0.\label{re5}\ee
By (\ref{re2}) and (\ref{re5}) there is   the smallest number $\rho_1\in[0,\rho_0)$ such that
\be\tau'>0\quad\mbox{and}\quad\var''(\rho)<0\qfq \rho\in(\rho_1,\rho_0).\label{re4}\ee
We claim $\rho_1=0.$ If $\rho_1>0,$ then by (\ref{re4}) and (\ref{re1})
$$\var'(\rho_1)>\var'(\rho_0)>\a_1(\rho_0,\var(\rho_0))\tau(\rho_0)\geq\a_1(\rho_1,\var(\rho_1))\tau(\rho_1),$$
since $\a_1(\rho_0,\var(\rho_0))\geq\a_1(\rho_1,\var(\rho_1))$ by $\rho_0\geq\rho_1$ and $\var(\rho_0)\geq\var(\rho_1).$
Using the same arguments as
for (\ref{re2}) and (\ref{re3}), we obtain
$$\tau'(\rho_1)>0,\quad\var''(\rho_1)<0,$$ reaching a contradiction.

It follows from (\ref{re4}) with $\rho_1=0$ that $\var'(\rho)>\var'(\rho_0)$ for $\rho\in(0,\rho_0)$ and
\be \var(\rho)<\var(\rho_0)-\var'(\rho_0)\rho_0+\var'(\rho_0)\rho\qfq\rho\in(0,\rho_0)\label{re6}.\ee
Moreover, (\ref{re1}) implies that
\be \frac{\rho_0\var'(\rho_0)}{\var(\rho_0)}>\a_1(\rho_0,\var(\rho_0))\frac{\rho_0}{f(\rho_0)}\frac{f\circ\var(\rho_0)}{\var(\rho_0)}
\geq\a_1(\rho_0,\var(\rho_0))\frac{b_0\circ\var(\rho_0)}{b_1(\rho_0)}\geq1,\label{re7}\ee since $\rho b_0(\rho)\leq f(\rho)\leq b_1(\rho)\rho$ for all $\rho\geq0.$

From (\ref{re6}) and (\ref{re7}) we obtain $\var(0)<0$ which is a contradiction again the assumption (\ref{re1}).

Similar arguments prove the left hand side of the inequalities (\ref{re}).  \hfill$\Box$\\

Actually, from the proof of Theorem \ref{t3.2} we have shown that

\begin{cor} \label{c3.1}Let $\mu_+(\max\{\lam,1\})\leq1.$ Let $(A 1),$ $(A 4),$ $(A 5),$ and $(A 6)$ hold.  Let  $\var\in C^1(0,1]$ be an equilibrium solution to  problem $(\ref{3.4}).$ Then
$$\var'(\rho)\leq \a_1(\rho,\var(\rho))\tau(\rho)\qfq\rho\in(0,1].$$ Furthermore, if there is $\rho_0\in(0,1]$ such that
$$\var'(\rho_0)< \a_0(\rho_0,\var(\rho_0))\tau(\rho_0),$$ then $\var(0)>0.$
\end{cor}

Let $\kappa\equiv0.$ Then $f=\rho,$ $b_0=b_1=1,$ and
$$\a_0(\rho,s)=\a_1(\rho,s)=1\qfq\rho,\,\,s\in[0,\infty).$$
It follows from (\ref{3.4}) and (\ref{re}) that a regular equilibrium solution $\var$ satisfies
\be\left\{\begin{array}{l}\var(1)=\lam,\\
\var'(\rho)=\dfrac{\var(\rho)}{\rho}\qfq\rho\in(0,1).\end{array}\right.\label{re8}\ee
Since  problem (\ref{re8}) has the unique solution $\var=\lam\rho,$ we have

\begin{cor}$(\cite{Si})$\label{c3.2}\,\,\,Let $\kappa\equiv0$ and let  $(A\,1),$ $(A\,4),$ $(A\,5),$ and $(A \,6)$ hold. Then there exists a unique regular equilibrium solution $\var=\lam\rho$ to problem
$(\ref{3.4})$ with $\var(1)=\lam.$
\end{cor}

Let $0<\mu_0\leq\mu_1<\infty$ be given by (\ref{pappa12}). It follows from (\ref{pappa12}) that
$$\frac{\mu_0}{\mu_1}\leq\frac{b_0(\rho)}{b_1(s)}\leq1\leq\frac{b_1(\rho)}{b_0(s)}\leq\frac{\mu_1}{\mu_0}\qfq\rho,\,\,s\in[0,\infty),$$
which also imply by (A 5) that
$$q_0^{-1}\Big(\frac{\mu_1}{\mu_0}\Big)\leq q_0^{-1}\Big(\frac{b_1(\rho)}{b_0(s)}\Big)\leq1\leq q_1^{-1}
\Big(\frac{b_0(\rho)}{b_1(s)}\Big)\leq q_1^{-1}\Big(\frac{\mu_0}{\mu_1}\Big)\qfq\rho,\,\,s\in[0,\infty).$$

It follows from Theorem \ref{3.2} and Proposition \ref{np1.1} that

\begin{cor} \label{c3.3} Let $\mu_+(\infty)\leq1$ and $\mu_-(\infty)<\infty.$ Let $(A \,1),$ $(A\,4),$ $(A \,5),$ and $(A \,6)$ hold.
Let $\var\in C^1(0,1]$ be a regular equilibrium solution to problem $(\ref{3.4})$ with $\var(1)=\lam.$ Then
\be  \eta_0\tau(\rho)\leq\var'(\rho)\leq\eta_1\tau(\rho)\qfq\rho\in(0,1],\label{re9}\ee where
$$\eta_0=\min\Big\{\,\frac{\mu_0}{\mu_1},\,q_0^{-1}\Big(\frac{\mu_1}{\mu_0}\Big)\,\Big\},\quad
\eta_1=\max\Big\{\,\frac{\mu_1}{\mu_0},\,q_1^{-1}\Big(\frac{\mu_0}{\mu_1}\Big)\,\Big\}. $$
\end{cor}

\begin{cor} \label{c3.4} Let all the assumptions in Corollary $\ref{c3.3}$ hold. If $\var\in C^1(0,1]$ is a regular equilibrium solution to  problem $(\ref{3.4})$ with $\var(1)=\lam,$ then
\be\lam\rho^{c_1}\leq\var(\rho)\leq\lam\rho^{c_0}\qfq \rho\in(0,1],\quad\lam>0,\label{3.24}\ee where
$$c_0=\eta_0\frac{\mu_0}{\mu_1},\quad c_1=\eta_1\frac{\mu_1}{\mu_0}.$$
\end{cor}

{\bf Proof}\,\,\,Using (\ref{pappa12}) and (\ref{re9}), we have
$$\frac{\var'}{\var}\leq\eta_1\frac{f\circ\var}{\var}\frac{\rho}{f}\frac{1}{\rho}\leq \frac{c_1}{\rho}\qfq\rho\in(0,1],$$ which yields the left hand side of the inequalities
(\ref{3.24}). A similar argument proves the right hand side of the inequalities
(\ref{3.24}). \hfill$\Box$

\begin{rem}
If $\kappa(s)=0$ for all $s\geq0,$ then $\mu_0=\mu_1=\eta_0=\eta_1=1,$ and $(\ref{3.24})$ means $\var(\rho)=\lam\rho.$
\end{rem}

\begin{thm} \label{t3.3}Let $\mu_+(\max\{\lam,1\})\leq1.$ Let $(A \,1),$  $(A \,4),$ $(A \,5),$ and $(A\,6)$  hold.
 Let  $\var\in C^1(0,1]$ be a regular equilibrium solution to problem $(\ref{3.4})$ with $\var(1)=\lam.$ Then

 $(i)$\,\,\, there are constants $\rho_0\in(0,1],$ $c_0>0,$ and $c_1>0$ such that
 \be c_0\rho\leq\var(\rho)\leq c_1\rho\quad\mbox{for all}\quad 0\leq\rho\leq\rho_0;\label{re3.26}\ee

$(ii)$\,\,\,the limit $\lim_{\rho\rw0+}\var(\rho)/\rho=\var'(0)$ exists;

 $(iii)$\,\,\, there are constants $\lam_0>0,$ $c_0>0,$ and $c_1>0$ such that
\be c_0\lam\leq\var'(0)\leq c_1\lam\quad\mbox{for all}\quad \lam\in(0,\lam_0].\label{re3.27*}\ee
\end{thm}

{\bf Proof}\,\,\,(i)\,\,\,First, we prove the right hand side of inequalities (\ref{re3.26}).

We fix $0<\rho_0\leq1$ small such that
\be e^{\mu_-\circ\var(\rho_0)}\mu_+\circ\var(\rho_0)<1,\label{re3.27}\ee   which is possible since $\var(0)=0.$

{\bf Step 1}\,\,\,Let
$$\a(\rho)=\frac{b_1(\rho)}{b_0\circ\var(\rho)}\qfq\rho\in(0,\rho_0].$$
We shall estimate
$$\int_\rho^{\rho_0}\frac{\a(s)}{\var(s)}\tau(s)ds\qfq\rho\in(0,\rho_0].$$

Since by Proposition \ref{p1.1}
$$b_0(\rho)=1-\sup_{0\leq s\leq\rho}\int_0^s\kappa fds\geq1-\int_0^\rho\kappa_+fds\geq 1-e^{\mu_-(\rho)}\mu_+(\rho)\qfq\rho\in(0,1],$$
it follows from (\ref{re3.27}) that
\be b_0\circ\var(\rho)\geq b_0\circ\var(\rho_0)>0\qfq 0\leq\rho\leq\rho_0.\label{re3.29}\ee

Moreover, we have
\beq &&b_1(\rho)-b_0\circ\var(\rho)=\sup_{0\leq s\leq\var(\rho)}\int_0^s\kappa fds-\inf_{0\leq s\leq\rho}\int_0^s\kappa fds\nonumber\\
&&\leq\int_0^{\var(\rho)}\kappa_+fds+\int_0^\rho\kappa_-fds\leq c(\rho_0,\var(\rho_0))[f\circ\var(\rho)+f(\rho)]
\label{re3.28}\eeq for $0\leq\rho\leq\rho_0,$ where
\be c(t,s)=\int_0^{\max\{t,s\}}|\kappa(\zeta)|d\zeta.\label{re39}\ee Using (\ref{re}), (\ref{re3.29}), and (\ref{re3.28}) we have
\beq &&\frac{\a(\rho)-1}{\var(\rho)}\tau(\rho)=\frac{\a(\rho)-1}{f(\rho)}\frac{f\circ\var(\rho)}{\var(\rho)}\leq\frac{ c(\rho_0,\var(\rho_0))}{b_0\circ\var(\rho_0)}(\tau+1)b_1\circ\var(\rho_0)\nonumber\\
&&
\leq \frac{ c(\rho_0,\var(\rho_0))b_1\circ\var(\rho_0)}{\a_0(\rho_0,\var(\rho_0))b_0\circ\var(\rho_0)}\var'(\rho)
+\frac{ c(\rho_0,\var(\rho_0))b_1\circ\var(\rho_0)}{b_0\circ\var(\rho_0)}\qfq0\leq\rho\leq\rho_0,\label{re3.30}\eeq since $f\circ\var\leq (b_1\circ\var)\var.$

In addition we have
\beq \frac{f\circ\var(\rho)}{\var(\rho)}-1&&=-\int_0^{\var(\rho)}\kappa fds+\frac{1}{\var(\rho)}\int_0^{\var(\rho)} s\kappa fds\nonumber\\
&&\leq 2c(\rho_0,\var(\rho_0))f\circ\var(\rho)\qfq0\leq\rho\leq\rho_0,\nonumber\eeq which implies by (\ref{re}) that
\be \Big[\frac{f\circ\var(\rho)}{\var(\rho)}-1\Big]\frac{1}{f(\rho)}\leq \frac{2c(\rho_0,\var(\rho_0))}{\a_0(\rho_0,\var(\rho_0))}\var'(\rho)\qfq0<\rho\leq\rho_0.\label{re3.31}\ee

Furthermore,
\be \frac{1}{f(\rho)}-\frac{1}{\rho}=\frac{1}{f(\rho)}\Big[\int_0^\rho\kappa(1-\frac{s}{\rho}) fds\Big]\leq 2c(\rho_0,\var(\rho_0))\qfq0<\rho\leq\rho_0.\label{re3.32}\ee

From (\ref{re3.30}), (\ref{re3.31}), and (\ref{re3.32}), we obtain
\beq \frac{\a(\rho)}{\var(\rho)}\tau(\rho)&&=\frac{\a(\rho)-1}{\var(\rho)}\tau(\rho)+\Big[\frac{f\circ\var(\rho)}{\var(\rho)}-1\Big]\frac{1}{f}+\Big(\frac{1}{f(\rho)}-\frac{1}{\rho}\Big)
+\frac{1}{\rho}\nonumber\\
&&\leq c(\rho_0,\var(\rho_0))\Big[\frac{b_1\circ\var(\rho_0)}{b_0\circ\var(\rho_0)}+2\Big]\Big[\frac{\var'(\rho)}{\a_0(\rho_0,\var(\rho_0))}+1\Big]
+\frac{1}{\rho}\nonumber\eeq for $0<\rho\leq\rho_0,$ which yields
\be \int_\rho^{\rho_0}\frac{\a(s)}{\var(s)}\tau(s)ds\leq\eta(\rho_0,\var(\rho_0))+\ln\frac{\rho_0}{\rho}\qfq0<\rho\leq\rho_0,\label{re3.33} \ee where
\be\eta(t,s)= c(t,s)\Big[\frac{b_1(s)}{b_0(s)}+2\Big]\Big[\frac{s}{\a_0(t,s)}+t\Big].\label{re3.35}\ee

{\bf Step 2}\,\,\,Let
$$\a*(\rho)=\frac{b_0(\rho)}{b_1\circ\var(\rho)}\qfq0\leq\rho\leq\rho_0.$$ By (\ref{re3.41})  we have
$$q_1^{-1}(\a*(\rho))-1\leq\frac{1}{c*(\rho_0,\var(\rho_0))}[1-\a*(\rho)]\qfq0\leq\rho\leq\rho_0,$$ where
\be c*(t,s)=\inf_{b_0(t)/b_1(s)\leq \zeta\leq1}|q_1'(q_1^{-1}(\zeta))|.\label{re40}\ee
Similar arguments as in {\bf Step 1} give the estimate
\be\int_\rho^{\rho_0}q_1^{-1}(\a*(s))\frac{f\circ\var(s)}{f(s)\var(s)}ds\leq \eta*(\rho_0,\var(\rho_0))+\ln\frac{\rho_0}{\rho}\qfq0<\rho\leq\rho_0,\label{re3.39}\ee
where
$$\eta*(t,s)=c(t,s)\Big[\frac{b_1(s)}{c*(t,s)}+2\Big]\Big[\frac{s}{\a_0(t,s)}+t\Big],$$ where $c(t,s)$ and $c*(t,s)$ are given by (\ref{re39}) and (\ref{re40}), respectively.

{\bf Step 3}\,\,\,Using (\ref{re}), (\ref{re3.33}), and (\ref{re3.39}) we obtain
\be \var(\rho)\geq e^{-\max\{\eta(\rho_0,\var(\rho_0)),\eta*(\rho_0,\var(\rho_0))\}}\frac{\var(\rho_0)}{\rho_0}\rho\qfq0\leq\rho\leq\rho_0.\label{re3.38*}\ee

Similar arguments prove the left hand side of inequalities (\ref{re3.26}).\\

(ii)\,\,\,From (\ref{re3.38*})  we have
\be \underline{\lim}_{\rho\rw0+}\frac{\var(\rho)}{\rho}\geq e^{-\max\{\eta(\rho_0,\var(\rho_0)),\eta*(\rho_0,\var(\rho_0))\}}
\frac{\var(\rho_0)}{\rho_0}\label{re3.40*}\ee for $\rho_0\in(0,1]$ such that (\ref{re3.27}) holds. It follows from (\ref{re3.40*}) that
$$\underline{\lim}_{\rho\rw0+}\frac{\var(\rho)}{\rho}\geq \overline{\lim}_{\rho\rw0+}\frac{\var(\rho)}{\rho}$$
since $\lim_{\rho\rw0+}\eta(\rho,\var(\rho))=\lim_{\rho\rw0+}\eta*(\rho,\var(\rho))=0,$ which is what we need.

(iii)\,\,\, Let $\lam_0>0$ be small such that
$$ e^{\mu_-(\lam_0)}\mu_+(\lam_0)<1.$$ Let $0<\lam\leq\lam_0.$ Then $\rho_0$ in (\ref{re3.27}) can be taken as $1.$
It follows from (\ref{re3.38*}) that
$$\var'(0)\geq e^{-\max\{\eta(1,\lam),\eta*(1,\lam)\}}\lam\qfq \lam\in(0,\lam_0].$$
Since $\eta(1,\lam)$ and $\eta*(1,\lam)$ are bounded on $[0,\lam],$ we have shown that the left hand side of inequalities (\ref{re3.27*}).
Similar arguments show that the right hand side of (\ref{re3.27*}) holds.\hfill$\Box$

\begin{rem}In $(\ref{re3.26})$ constants $c_0$ and $c_0$ may be depend on the number $\lam.$

\end{rem}

\begin{pro}\label{p3.2}Let $\mu_+(\max\{\lam,1\})\leq1.$ Let $(A \,1),$ $(A\,2),$ $(A \,4),$ $(A \,5),$ $(A\,6),$ $(A\,7)$ and $(A\,8)$  hold.
 Let  $\var\in C^1(0,1]$ be an equilibrium solution to problem $(\ref{3.4})$ with $\var(1)=\lam.$
 Then

$(i)$\,\,\, $\sup_{0<\rho\leq1}\var'(\rho)<\infty;$

$(ii)$\,\,\, the limit $\lim_{\rho\rw0+}T(\rho)$ exists and is finite.
\end{pro}

{\bf Proof}\,\,\, If $\var$ is  regular, then (i) and (ii) follow from Theorems \ref{t3.3} and \ref{t3.2}.

Next, we assume that $\var(0)>0.$ Using Corollary \ref{c3.1}, (A 7), (A 8), and (\ref{3.3*}) we obtain
$$ |T'(\rho)|\leq c\tau^{1-n}[1+\tau^{\a}+(\a_1\tau)^{\a}]\leq c\tau^{1+\a-n}\leq \dfrac{c}{\rho^{1+\a-n}}\qfq\rho\in(0,1],$$ which yields
\be |T(\rho)|\leq|T(1)|+c\int_\rho^1\frac{1}{\rho^{1+\a-n}}d\rho\leq |T(1)|+c\qfq\rho\in(0,1]\label{ca3.41}\ee since $1+\a-n<1.$

To prove (i) we suppose for a contradiction that there is a sequence $\{\,\rho_j\,\}\subset(0,1]$ such that $\rho_j\rw0$ and
$$j\leq\var'(\rho_j)\leq \a_1(\rho_j,\var(\rho_j))\tau(\rho_j)\qfq j\geq1,$$ where the right hand side of the above inequalities is from Corollary \ref{c3.1}.
It follows from (A 6) and (A 7) that
$$\tau^{1-n}(\rho_j)|\phi'(\var'(\rho_j))|\leq c\tau^{1-n}(\rho_j)[1+\var'^{\a-1}(\rho_j)]\leq c[\tau^{1-n}(\rho_j)+\tau^{\a-n}(\rho_j)],$$ and then
$\lim_{j\rw\infty}\tau^{1-n}(\rho_j)\phi'(\var'(\rho_j))=0$ is true. By (A 2) we have
$$\lim_{j\rw\infty}T(\rho_j)=\lim_{j\rw\infty}h'\Big(\var'(\rho_j)\tau^{n-1}(\rho_j)\Big)=+\infty,$$ contradicting with (\ref{ca3.41}).

By (i) and $\var(0)>0$ there is some $\rho_0\in(0,1]$ such that
$$\var'(\rho)<\tau(\rho)\qfq\rho\in(0,\rho_0].$$  Then $T'(\rho)\geq0$ for $\rho\in(0,\rho_0]$ by (\ref{3.3*}), which implies that (ii) is true. \hfill$\Box$\\

{\bf Cavitating Equilibrium Solutions}\,\,\,An equilibrium solution $\var\in C^1(0,1]$ to problem (\ref{3.4}) is said to be {\it cavitating} if

(i)\,\,\,$\var(0)>0,$ and

(ii)\,\,\,$\lim_{\rho\rw0+}T(\rho)=0.$

 Ball \cite{Ball1} has shown
that when $\lam>0$ is small there is no cavitating equilibrium solution in the case of the Euclidean space. We  present some  similar results.

\begin{pro}\label{p3.3} Let $\mu_+(\max\{\lam,1\})\leq1.$  Let $(A \,1),$ $(A\,2),$ $(A \,4),$ $(A \,5),$ $(A\,6),$  and $(A\,7)$  hold.
Suppose that $(A\,8)$ holds with $\b=0.$
If $\var\in C^1(0,1]$ is a cavitating equilibrium solution to  problem $(\ref{3.4}),$
then
\be \lim_{\rho\rw0+}\var'\tau^{n-1}=\varpi,\label{ree3.42}\ee where $\varpi>0$ is given by $h'(\varpi)=0.$
\end{pro}

{\bf Proof}\,\,\,By Corollary \ref{c3.1},  (A 7), and (A 8) with $\b=0$, we have
$\lim_{\rho\rw0+}\tau^{n-1}\phi'(\var')=0.$ Then $0=\lim_{\rho\rw0+}T(\rho)=\lim_{\rho\rw0+}h'(\var'\tau^{n-1}).$
(\ref{ree3.42}) follows from (A 1). \hfill$\Box$

\begin{thm}\label{t3.4}
Let $\mu_+(1)\leq1$ and let there be $\varepsilon\in(0,1]$ such that
\be \kappa(s)=0\qfq s\in(0,\varepsilon).\label{3.25}\ee
Let $(A \,1),$ $(A\,2)$ and $(A \,4)$ hold.
Then there is no cavitating equilibrium solution when $\lam>0$ is small.
\end{thm}

{\bf Proof}\,\,\,It follows from (\ref{3.25}) that
\be f'(\rho)=1,\quad f(\rho)=\rho\qfq\rho\in[0,\varepsilon].\label{ca44}\ee Then
the left hand side of the equation in (\ref{3.4}) equals
\beq&&(n-1)[\phi'(\tau)
-\phi'(\var')-(\var'-\tau)h''(\var'\tau^{n-1})\var'\tau^{2n-3}]\label{ca45}\eeq
for $\rho\in(0,\varepsilon]$ when $\var(1)=\lam<\varepsilon,$ where
$$\tau(\rho)=\frac{\var(\rho)}{\rho}\qfq\rho\in(0,\varepsilon].$$

By (A 2) and (A 4), without loss of generality we assume that $\lam_0>0$ is small such that
\be \Phi_1(\frac{\lam_0}{\varepsilon},\cdots,\frac{\lam_0}{\varepsilon})=\phi'(\frac{\lam_0}{\varepsilon})+h'\Big((\frac{\lam_0}{\varepsilon})^n\Big)<0.\label{ca46}\ee

Let $\var\in C^1(0,1]$ be a cavitating equilibrium solution with $\var(1)=\lam<\lam_0$ and we derive a contradiction below.
We claim
\be \var'(\rho)<\tau(\rho)\quad\mbox{for all}\quad\rho\in(0,\varepsilon].\label{ca47}\ee In fact, if there were some $\rho_0\in(0,\varepsilon]$ such that
$\var'(\rho_0)=\tau(\rho_0),$ it is easy to check by (\ref{ca45}) and the uniqueness theorem that
$$\var(\rho)=\var'(\rho_0)\rho\quad\mbox{for all}\quad \rho\in(0,\rho_0],$$ which contradicts $\var(0)>0.$ Next, by (\ref{ca47}), (\ref{3.3*}) and (\ref{ca46}),
we obtain
\beq T(\rho)&&<T(\varepsilon)=\tau^{1-n}(\varepsilon)\Phi_1(\var'(\varepsilon),\tau(\varepsilon),\cdots,\tau(\varepsilon))\nonumber\\
&&<\tau^{1-n}(\varepsilon)\Phi_1(\tau(\varepsilon),\tau(\varepsilon),\cdots,\tau(\varepsilon))<0\qfq\rho\in(0,\varepsilon),\nonumber\eeq that contradicts $T(0)=0.$
\hfill$\Box$

\begin{thm}\label{t3.5} Let $\mu_+(1)\leq1$ and let there be $\varepsilon\in(0,1]$ such that
\be \kappa(s)<0\qfq s\in(0,\varepsilon).\label{3.25*}\ee
Let $(A \,1),$ $(A\,2),$ $(A \,4),$ $(A \,5),$ $(A\,6),$ $(A\,7),$  and $(A\,8)$ hold. Further suppose that

$(i)$\,\,\,$\phi'(1)+h'(1)=0;$

$(ii)$\,\,\,$\phi'(s)-h''(s^n)s^{2n-1}<0$ for $s>1.$

Then there is no cavitating equilibrium solution when $0<\lam<\varepsilon.$
\end{thm}

{\bf Proof}\,\,\,We suppose for a contradiction that $\var\in C^1(0,1]$ is a cavitating equilibrium solution with $\var(1)=\lam<\varepsilon.$ By Proposition \ref{p3.2} there is the largest number
$\rho_0\in(0,1]$ such that
\be \var'(\rho)<\tau(\rho)\quad\mbox{for all}\quad\rho\in(0,\rho_0).\label{ca3.44}\ee

Let $\rho_0\geq\varepsilon.$ Then
 $$\tau(\rho_0)\leq\frac{f(\lam)}{f(\varepsilon)}<1,$$  and therefore there is some $\rho_1\in(0,\rho_0)$ such that $\tau(\rho_1)=1.$ The relations (\ref{3.3*}) and (\ref{ca3.44}) implies that
\be T(\rho)\leq T(\rho_1)=\tau^{1-n}(\rho_1)\Phi_1(\var'(\rho_1),1,\cdots,1)<\tau^{1-n}(\rho_1)\Phi_1(1,\cdots,1)=0\label{ca3.46}\ee for all $\rho\in(0,\rho_1),$ which contradicts with $T(0)=0.$

Let $0<\rho_0<\varepsilon.$ By (\ref{ca3.44}) we must have
\be \var'(\rho_0)=\tau(\rho_0).\label{ca3.45}\ee If $\tau(\rho_0)\leq1,$ we get a contradiction as in (\ref{ca3.46}). We assume that $\tau(\rho_0)>1.$ Then
\be \rho_0<\var(\rho_0)<\lam<\varepsilon.\label{ca3.47}\ee
By (\ref{ca3.47}) and (\ref{3.25*}), we obtain
\be f'\circ\var(\rho_0)=f'(\rho_0)-\int_{\rho_0}^{\var(\rho_0)}\kappa(s)f(s)ds>f'(\rho_0).\label{ca53}\ee Thus, it follows from (\ref{3.4}), (\ref{ca3.45}), (\ref{ca53}) and (ii) that
\beq &&\mbox{the left hand side of the equation (\ref{3.4})}\nonumber\\
&&=(n-1)(f'\circ\var-f')[\phi'(\tau)-h''(\tau^n)\tau^{2(n-1)}]<0\quad\mbox{at}\quad \rho=\rho_0,\nonumber\eeq that is,
\be\var''(\rho_0)<0. \label{ca54}\ee In addition, we have, by (\ref{ca53}) and (\ref{ca3.45}),
\be f(\rho_0)\tau'(\rho_0)=[f'(\var(\rho_0)-f(\rho_0)]\tau(\rho_0)>0.\label{ca55}\ee Thus, by (\ref{ca54}) and (\ref{ca55}), we obtain
$$\var'(\rho)>\var'(\rho_0)=\tau(\rho_0)>\tau(\rho)\quad\mbox{for}\quad \rho<\rho_0\quad\mbox{and near}\quad \rho_0,$$ contradicting the definition of $\rho_0.$ \hfill$\Box$

\begin{rem} If $\kappa(o)<0,$ assumption $(\ref{3.25*})$ is true.
\end{rem}

\begin{rem} Assumption $(i)$ means that $\Phi_i(1,\cdots,1)=0$  for $1\leq i\leq n$ so that the undeformed configuration
is a natural state. In addition,  if
$$\phi(v)=\frac{1}{v^\b}+ v^\a,\quad h(v)=\frac{2\a}{v}+\frac{\a+\b}{2} v^2,\quad \b>0,\quad\a>1,$$ then $(i)$ and $(ii)$ hold.
\end{rem}

\subsection{Energy Minimizers}
\hskip\parindent We seek to minimize
\be I(\var)=\int_0^1f^{n-1}(\rho)\Phi(\rho)d\rho\label{4.1}\ee among radial deformations $\var$ such that $\var(0)\geq0$ and $\var(\rho)$ is increasing.
For $\lam>0$ given, let
$$\A_\lam=\Big\{\,\var\in W^{1,1}(0,1)\,\Big|\,\var(0)\geq0,\,\var'>0\,\,\mbox{a.e., and}\,\,I(\var)<\infty\,\Big\}.$$ Consider the minimization problem
\be\inf_{\var\in\A_\lam}I(\var).\label{4.2}\ee

Let
$$u(p)=\si\circ\var(\rho),$$ where $\si$ is given by (\ref{si}) and $p=\si(\rho).$ Then
$$\var(\rho)=\si^{-1}\circ u(p),\quad \rho=\si^{-1}(p),$$
$$ \var'(\rho)=\frac{u'(p)}{\si'\circ\var(\rho)}p'(\rho)=\frac{f^{n-1}\circ\si^{-1}(p)}{f^{n-1}\circ\si^{-1}(u)}u'(p),\quad \tau(\rho)=\frac{f\circ\si^{-1}(u)}{f\circ\si^{-1}(p)},$$
and
$$u'(p)=\si'\circ\var(\rho)\var'(\rho)\frac{\pl\rho}{\pl p}=\var'(\rho)\tau^{n-1}(\rho).$$
Let
$$\wp(p,u,q)=\phi\Big(\frac{f^{n-1}\circ\si^{-1}(p)}{f^{n-1}\circ\si^{-1}(u)}q\Big)+(n-1)\phi\Big(\frac{f\circ\si^{-1}(u)}{f\circ\si^{-1}(p)}\Big)+h(q),$$
$$J(u)=\int_0^{\si(1)}\wp(p,u,u')dp.$$ Since
$$I(\var)=J(u),$$ the minimization problem (\ref{4.2}) is now equivalent to minimizing  $J(u)$ on the set
$$\Big\{\,u\in W^{1,1}(0,\si(1))\,\Big|\,u(\si(1))=\si(\lam),\,u(0)\geq0,\,u'(p)>0\,\mbox{a.e., and}\,J(u)<\infty\,\Big\}.$$

A similar argument as in \cite{Ball1} yields

\begin{thm}\label{t4.1} Suppose that $(A\,1),$ $(A\,2),$ $(A\,4),$ and $(A\,11)$ hold.
$I$ attains an absolute minimum on $\A_\lam.$
\end{thm}

Following the proof of  Theorem 7.3 in \cite{Ball1},  we obtain the following.

\begin{thm}\label{t4.2}Let  $(A\,1),$ $(A\,2),$ $(A\,4),$ $(A\,7),$ and $(A\,9)$ hold. Let $\var$ be a minimizer of $I$ on $\A_\lam.$ Then $\var\in C^1(0,1],$
$\var'(\rho)>0$ for all $\rho\in(0,1],$ $f^{n-1}\Phi_1\in C^1(0,1]$ and $(\ref{3.4})$ holds for all $\rho\in(0,1].$ If $\var(0)>0,$ then $f^{n-2}\Phi_2\in L^1(0,1)$ and
\be\lim_{\rho\rw0+}T(\rho)=0,\label{4.3}\ee where $T$ is the radial component of the Cauchy stress, given by $(\ref{3.2*}).$
\end{thm}

Next, we have the following.

\begin{thm}Let $\mu_+(\max\{\lam,1\})\leq1.$ Let $(A \,1),$ $(A\,2),$ $(A \,4),$ $(A \,5),$ $(A\,6),$ $(A\,7)$ and $(A\,8)$  hold.
 Let  $\var\in C^1(0,1]$ be an equilibrium solution to  problem $(\ref{3.4})$ with $\var(1)=\lam.$ Then
 $$I(\var)<\infty.$$
\end{thm}

{\bf Proof}\,\,\,By using (\ref{3.2}) and (\ref{3.3*}), we have
\beq\var'\Phi_1+(n-1)\tau\Phi_2&&=(n-1)[\tau\phi'(\tau)-\var'\phi'(\var')]+n\var'\tau^{n-1}T\nonumber\\
&&=\frac{f\circ\var}{f'\circ\var}\tau^{n-1}T'+n\var'\tau^{n-1}T\qfq\rho\in(0,1).\nonumber\eeq
Thus,
\be f^{n-1}[\var'\Phi_1+(n-1)\tau\Phi_2]=\frac{1}{f'\circ\var}\Big(f^n\circ\var T\Big)'\qfq\rho\in(0,1].\label{ca4.4}\ee
Now, using (\ref{3.7*}) and (\ref{ca4.4}), we obtain, by integration by parts,
\beq &&n\int_\rho^1f'f^{n-1}\Phi ds+f^n(\rho)\Phi(\rho)=\Big\{f^n[\Phi-(\var'-\tau)\Phi_1]+\frac{f'-f'\circ\var}{f'\circ\var}f^n\circ\var T\Big\}\Big|_{\rho=1}\nonumber\\
&&\quad+\frac{f'\circ\var(\rho)-f'(\rho)}{f'\circ\var(\rho)}f^n\circ\var(\rho)T(\rho)+f^n(\rho)[\var'(\rho)-\tau(\rho)]\Phi_1(\rho)\nonumber\\
&&\quad+\int_\rho^1\frac{\kappa\circ\var\var'f'-\kappa f f'\circ\var}{f'^2\circ\var}f^n\circ\var Tds\qfq\rho\in(0,1].\label{ca4.5}\eeq
Next, by Proposition \ref{p3.2}, the term
\beq f^n(\rho)[\var'(\rho)-\tau(\rho)]\Phi_1(\rho)&&=f^{n-1}\circ\var(\rho)[f(\rho)\var'(\rho)-f\circ\var(\rho)]T(\rho)\nonumber\eeq\
converges as $\rho$ goes to $0+$ and the limit is finite.
Moreover, by Proposition \ref{p1.1},
$$\mu_0(1)\leq f'(\rho)\leq e^{\mu_1(1)},\quad\mu_0(\lam)\leq f'\circ\var(\rho)\leq e^{\mu_1(\lam)}\qfq\rho\in(0,1].$$
Thus, the second term and the last term in the right hand side of (\ref{ca4.5}) converge as $\rho$ goes to $0+$ and their limits are finite by Proposition \ref{p3.2}.
The proof is complete.  \hfill$\Box$

\subsection{Cavitation}
\hskip\parindent
We now derive a cavitating theorem.

\begin{thm}\label{t4.4} Let $\mu_+(\infty)\leq1$ and $\mu_-(\infty)<\infty.$
Suppose that $(A\,1)$-$(A\,9)$ hold. For any $\lam$ sufficiently large a minimizer of $I$ on $\A_\lam$ is a cavitating equilibrium solution.
\end{thm}

{\bf Proof}\,\,\,By Theorem \ref{t4.2} a minimizer of $I$ on $\A_\lam$ is an equilibrium solution. We suppose for a contradiction that  $\var\in C^1(0,1]$ is
a regular equilibrium solution with $\var(1)=\lam$ where $\lam$ is sufficiently large. For $\varepsilon\in(0,1),$ let
$$\var_\varepsilon(\rho)=\var\Big\{\si^{-1}[(1-\varepsilon)\si(\rho)+\varepsilon\si(1)]\Big\}.$$ Then
$$\var_\varepsilon(1)=\lam,\quad\var_\varepsilon(0)=\var(\si^{-1}(\varepsilon\si(1)))>0. $$
Thus, $\var_\varepsilon\in\A_\lam$ for all $\varepsilon\in(0,1).$ We will show that for $\lam$ sufficiently large and $\varepsilon\in(0,1)$ small
\be I(\var_\varepsilon)<I(\var),\ee contradicting that $\var$ is a minimizer of $I$ on $\A_\lam.$

{\bf Step 1}\,\,\,Let
$$p=\si^{-1}[(1-\varepsilon)\si(\rho)+\varepsilon\si(1)],\quad\tau_\varepsilon=\frac{f\circ\var_\varepsilon(\rho)}{f(\rho)},\qfq\rho\in(0,1].$$ Simple computations yield
\be p'=(1-\varepsilon)\frac{f^{n-1}}{f^{n-1}\circ p},\quad\var_\varepsilon'(\rho)=(1-\varepsilon)\var'\circ p(\rho)\frac{f^{n-1}(\rho)}{f^{n-1}\circ p(\rho)},\label{ca4.17}\ee
\be\var_\varepsilon'(1)=(1-\varepsilon)\var'(1),\quad\tau_\varepsilon=\tau\circ p\frac{f\circ p}{f},\quad \var_\varepsilon'\tau^{n-1}_\varepsilon=(1-\varepsilon)\var'\circ p\tau^{n-1}\circ p.\label{ca4.18}\ee
In addition, (A 1) gives
$$h(t)\geq h(s)+(t-s)h'(s)\quad\mbox{for all}\quad t,\,\,s\in(0,\infty).$$ By using the formulas above, we obtain
\beq &&(1-\varepsilon)\int_0^1h(\var'_\varepsilon\tau^{n-1}_\varepsilon)f^{n-1}d\rho=\int_{p(0)}^1h[(1-\varepsilon)\var'\tau^{n-1}]f^{n-1}d\rho\nonumber\\
&&\leq \int_{p(0)}^1h(\var'\tau^{n-1})f^{n-1}d\rho-\varepsilon\int_{p(0)}^1\var'\tau^{n-1}h'[(1-\varepsilon)\var'\tau^{n-1}]f^{n-1}d\rho,\nonumber\eeq
which yields
\beq \int_{0}^1h(\var'\tau^{n-1})f^{n-1}d\rho&&\geq \varepsilon\int_{0}^{p(0)}h(\var'\tau^{n-1})f^{n-1}d\rho+(1-\varepsilon)\int_0^1h(\var'_\varepsilon\tau^{n-1}_\varepsilon)f^{n-1}d\rho\nonumber\\
&&\quad
+\varepsilon(1-\eta)\int_{p(0)}^1\var'\tau^{n-1}h'[(1-\varepsilon)\var'\tau^{n-1}]f^{n-1}d\rho\nonumber\\
&&\quad+\varepsilon\eta\int_{p(0)}^1\var'\tau^{n-1}h'[(1-\varepsilon)\var'\tau^{n-1}]f^{n-1}d\rho,\label{ca4.7}\eeq
for any $\eta>0$ small, since $\varepsilon<1.$

{\bf Step 2}\,\,\,By (A 3), there are constants $\varpi$ and $N(\varpi)>0$ such that
\be \varpi>1,\quad xh'(x)\geq\varpi h(x),\quad\mbox{for all}\quad x\geq N(\varpi).\label{ca4.8}\ee
Since $\theta$ is continuous, we take $N(\varepsilon)\geq N(\varpi)$ such that
\be h[(1-\varepsilon)x]\geq[\theta(1-\varepsilon)-\varepsilon]h(x)\quad\mbox{for all}\quad x\geq N(\varepsilon).\ee Next, for $\varpi>1$ being given in (\ref{ca4.8}),
we fix $\eta>0$ and $\varepsilon>0$ small such that
\be \frac{(1-\eta)\varpi}{1-\varepsilon}[\theta(1-\varepsilon)-\varepsilon]\geq1,\label{ca4.10} \ee that is possible because $\lim_{\varepsilon\rw0}\theta(1-\varepsilon)=1.$

By using (\ref{ca4.8})-(\ref{ca4.10}), we obtain that, if $\var'\tau^{n-1}\geq N(\varepsilon)/(1-\varepsilon)$ for all $\rho\in[p(0),1],$ then
\beq &&(1-\eta)\int_{p(0)}^1\var'\tau^{n-1}h'[(1-\varepsilon)\var'\tau^{n-1}]f^{n-1}d\rho\geq\frac{\varpi}{1-\varepsilon}\int_{p(0)}^1h[(1-\varepsilon)\var'\tau^{n-1}]f^{n-1}d\rho
\nonumber\\
&&\geq\frac{(1-\eta)\varpi}{1-\varepsilon}[\theta(1-\varepsilon)-\varepsilon]\int_{p(0)}^1h(\var'\tau^{n-1})f^{n-1}d\rho\geq \int_{p(0)}^1h(\var'\tau^{n-1})f^{n-1}d\rho.\label{ca4.11}\eeq
Inserting (\ref{ca4.11}) into (\ref{ca4.7}) yields
\beq \int_0^1h(\var'\tau^{n-1})f^{n-1}d\rho&&\geq \int_0^1h(\var'_\varepsilon\tau_\varepsilon^{n-1})d\rho\nonumber\\
&&\quad+\frac{\varepsilon\eta}{1-\varepsilon}\int_{p(0)}^1\var'\tau^{n-1}h'[(1-\varepsilon)\var'\tau^{n-1}]f^{n-1}d\rho,\label{ca4.12}\eeq if the following condition holds
\be \var'\tau^{n-1}\geq\frac{N(\varepsilon)}{1-\varepsilon}\quad\mbox{for all}\quad \rho\in[p(0),1].\ee

{\bf Step 3}\,\,\,Since $\var$ is regular, by Corollaries \ref{c3.3}, \ref{c3.4} and Proposition \ref{np1.1}, we obtain
\be\lam\frac{\mu_0}{\mu_1}\rho^{c_1-1}\leq\tau(\rho)\leq\lam\frac{\mu_1}{\mu_0}\rho^{c_0-1},\quad
\lam\frac{\eta_0\mu_0}{\mu_1}\rho^{c_1-1}\leq\var'(\rho)\leq\lam\frac{\eta_1\mu_1}{\mu_0}\rho^{c_0-1},\label{ca4.14}\ee for all $\rho\in(0,1].$

By (A 8) and (\ref{ca4.14}), we have
\beq \phi(\var')+(n-1)\phi(\tau)&&\leq c(p(0))(1+\lam^\a),\qfq\lam\geq1,\,\,\rho\in[p(0),1],\eeq where constant $c(p(0))>0$ depends on $p(0)$ but is independent of $\lam\geq1.$
Next, by (\ref{ca4.17}), (\ref{ca4.18}) and (\ref{ca4.14}), we obtain
\be \lam(1-\varepsilon)\frac{\eta_0\mu_0}{\mu_1}p^{c_1-1}(0)\frac{f^{n-1}(\rho)}{f^{n-1}(1)}\leq\var_\varepsilon'(\rho)
\leq\lam(1-\varepsilon)\frac{\eta_1\mu_1}{\mu_0}p^{c_0-1}(0)\frac{f^{n-1}(\rho)}{f^{n-1}\circ p(0)}\label{ca4.19}\ee for all $\rho\in(0,1],$ and
\be\lam p^{c_1-1}(0)\frac{\mu_0}{\mu_1}\frac{f\circ p(0)}{f(\rho)}\leq\tau_\varepsilon(\rho)\leq\lam\frac{\mu_1}{\mu_0}p^{c_0-1}(0)\frac{f(1)}{f(\rho)}
\quad\mbox{for all}\quad\rho\in(0,1]\label{ca4.20} \ee and
\be \var'(\rho)\tau^{n-1}(\rho)\geq\lam^n \hat{c}(p(0))\qfq\rho\in[p(0),1].\ee It follows from (A 8), (\ref{ca4.19}) and (\ref{ca4.20}) that
\be \phi(\var'_\varepsilon)+(n-1)\phi(\tau_\varepsilon)\leq c(p(0))\Big(1+\lam^\a[1+\frac{1}{f^\a(\rho)}]+\frac{1}{f^{\b(n-1)}(\rho)}+\frac{1}{f^\b(\rho)}\Big)\label{ca4.21}\ee for
all $\lam\geq1$ and $\rho\in(0,1].$

By using (A 4), (\ref{ca4.12}),  (A 8) and (\ref{ca4.19})-(\ref{ca4.21}),  we obtain
\beq I(\var)&&\geq\int_{p(0)}^1[\phi(\var')+(n-1)\phi(\tau)]f^{n-1}d\rho+\int_0^1h(\var'\tau^{n-1})f^{n-1}d\rho\nonumber\\
&&\geq I(\var_\varepsilon)+\int_{p(0)}^1[\phi(\var')+(n-1)\phi(\tau)]f^{n-1}d\rho\nonumber\\
&&\quad-\int_{0}^1[\phi(\var'_\varepsilon)+(n-1)\phi(\tau_\varepsilon)]f^{n-1}d\rho\nonumber\\
&&\quad+\frac{\varepsilon\eta}{1-\varepsilon}\int_{p(0)}^1\var'\tau^{n-1}h'[(1-\varepsilon)\var'\tau^{n-1}]f^{n-1}d\rho\nonumber\\
&&\geq I(\var_\varepsilon)+\frac{\varepsilon\eta}{1-\varepsilon}\lam^n\hat{c}(p(0))h'[(1-\varepsilon)\lam^n\hat{c}(p(0))]
\int_{p(0)}^1f^{n-1}(\rho)d\rho\nonumber\\
&&\quad-c(p(0))(1+\lam^\a),\eeq when
$$\lam\geq\max\Big\{\Big(\frac{N(\varepsilon)}{\hat{c}(p(0))(1-\varepsilon)}\Big)^{1/n},\,\,\,1\Big\}.$$

The proof is complete. \hfill$\Box$

\setcounter{equation}{0}
\section{Cavitation for Membrane Shells}
\def\theequation{5.\arabic{equation}}
\hskip\parindent The nonlinear shell membrane energy is obtained in \cite{LeRa} by the $\Ga-$ limit of the sequence of three-dimensional energies.
Here we will show that such membrane enerygies may take a form of (\ref{4}) where $M$ is a surface in $\R^3$ and $g$ is the induced metric of $M$ from
$\R^3.$

We consider a homogeneous elasitic material with stored-energy function $\hat{W}:$ $M_+^{3\times3}\rw\R.$ Suppose that $\hat{M}$ is frame-indifferent and isotropic, that is,
\be \hat{W}(F)=\hat{W}(QFR)\qfq F\in M_+^{3\times3},\quad Q,\,\,R\in \SOIII.\label{6.1}\ee

We assume that a middle surface $S$ is a bounded, connected open set of a $C^2$ surface in $R^3$  and let $N$ be the normal field of $S.$ For $h>0$ given, we consider the set
$\Om_h$ defined by
$$\Om_h=\{\,x+sN(x)\,|\,x\in S,\,\,|s|<h\,\}.$$ This set is the reference configuration of a shell with thickness $2h.$  Let $\u:$ $\Om_h\rw\R^3$ be a deformation of the shell. Then the stored energy is
\beq E_h(\u)&&=\int_{\Om_h}\hat{W}(\na\u(x))dx=\int_{-h}^h\int_S\hat{W}(\na\u(x+sN(x)))A(x,s)dgds,\nonumber\eeq where $g$ is the induced metric on $S$ from $\R^3$ and
$$A(x,s)=1+sH+s^2\kappa,$$ $H/2$ is the mean curvature, and $\kappa$ is the Gaussian curvature in $S.$

 For a deformation $\u\in L^p(\Om_1,\R^3),$ we define $\u_h\in L^p(\Om_h,\R^3)$ by
\be\u_h(x+shN(x))=\u(x+sN(x))\qfq x+shN(x)\in\Om_h.\label{6.3}\ee We define the rescaled energies by
\beq I_h(\u)=\frac{1}{h}E_h(\u_h)&&=\int_{-1}^1\int_S\hat{W}(\na\u_h(x+hsN(x)))A(x,sh)dgds,\nonumber\eeq  for $\u\in L^p(\Om_1,\R^3)$ and   $h>0$ small.

Let $x\in S$ be fixed. Let $e_1,$ $e_2$ be an orthonormal basis of $S_x$ such that $e_1,$ $e_2,$ $e_3$ is an orthonormal basis of $\R^3$ with positive orientation,
where $e_3=N(x),$ to satisfy
\be \hat{D}_{e_i}N=\kappa_ie_i\quad\mbox{at}\quad x\qfq i=1,\,\,2,\ee where $\hat{D}$ is the connection of the Euclidean space
$\R^3$ and $\kappa_i$ are the eigenvalues of the second fundamental form of $S$ at $x.$ Thus, $H(x)=\kappa_1+\kappa_2$
and $\kappa(x)=\kappa_1\kappa_2.$

It follows from (\ref{6.3}) that
\be(1+sh\kappa_i)\u_{h*}e_i=(1+s\kappa_i)\u_*e_i,\quad h\u_{h*}e_3=\u_*e_3.\label{n6.4}\ee
By (\ref{6.1}) and (\ref{n6.4}), we obtain
\beq \hat{W}(\na\u_h(x+hsN(x)))&&=\hat{W}(\,\u_{h*}e_1\,|\,\u_{h*}e_2\,|\u_{h*}e_3\,)\nonumber\\
&&=\hat{W}\Big(\,\frac{1+s\kappa_1}{1+sh\kappa_1}\u_{*}e_1\,\Big|\,\frac{1+s\kappa_2}{1+sh\kappa_2}\u_{*}e_2\,\Big|\,\frac{1}{h}\u_{*}e_3\Big).\eeq
In particular, if we let
$$\var^h(x+sN(x))=\var(x)+shw(x)\qfq \var,\,\,w\in W^{1,p}(S,\R^3),$$ then
$$\hat{W}(\na\var^h_h(x+shN))=\hat{W}\Big(\frac{\var_*e_1+shw_*e_1}{1+sh\kappa_1}\,\Big|\,\frac{\var_*e_2+shw_*e_2}{1+sh\kappa_2}\,\Big|\,w\Big).$$

Let
$$\hat{W}_0(F_1,F_2)=\min_{z\in\R^3}\hat{W}(F_1|F_2|z)\qfq F_1,\,\,F_2\in\R^3.$$
Let $\Q \hat{W}_0=\sup\{\,Z:\,R^3\times\R^3\rw\R,\,Z\,\mbox{quasiconvex},\,Z\leq \hat{W}_0\,\}$ be the
quasiconvex envelope of $\hat{W}_0.$

We introduce the space
\be V_M=\{\,\var\in W^{1,p}(\Om_1,\R^3)\,|\,\var_*N=0\,\mbox{for}\, x\in S\,\},\ee
for which we call the space of membrane displacements.

By similar arguments as in \cite{LeRa}, we may have
\begin{thm} \label{nt5.1}$(\cite{LeRa})$ Let all the assumptions on the function $\hat{W}$ in $\cite{LeRa}$ hold. Further suppose $(\ref{6.1})$ is true. Then the sequence $I_h$ $\Ga$-converges for the strong topology of $L^p(\Om_1,\R^3)$ when $h\rw0.$ For $\var\in L^p(\Om_1,R^3)$, it's
$\Ga$-limit is given by
\be I_0(\var)=\left\{\begin{array}{l}2\int_S\Q \hat{W}_0(\var_*e_1,\var_*e_2)dg\quad\mbox{if}\quad\var\in V_M,\\
+\infty\quad\mbox{otherwise.}\end{array}\right.\label{6.7}\ee
\end{thm}

Next, we shall reformulate   the $\Ga$-limit energy (\ref{6.7}) to relate it to the energy formula (\ref{4}).

Let us introduce a function $W_0:$ $M_+^{2\times2}\rw\R$ by
$$W_0(F)=\min_{z\in\R^3}\hat{W}\left(\begin{array}{cc}F&Z_1\\
0&z_3\end{array}\right),$$where $z=(z_1,z_2,z_3)$ and $Z_1=(z_1,z_2)^T.$ It is easy to check that
\be W_0(QFR)=W_0(F)\qfq F\in M_+^{2\times2},\quad Q,\,\,R\in\SOII.\label{6.5}\ee

Let $ W=\sup\{\,Z:\,M^{2\times2}\rw\R,\,Z\,\mbox{quasiconvex},\,Z\leq W_0\,\}$ be the
quasiconvex envelope of $W_0.$
From Dacorogna's representation formula for the quasiconvex envelope of $W_0,$ we obtain
$$W(F)=\frac{1}{\pi}\inf_{\X\in W_0^{1,\infty}(D,\R^2)}\int_DW_0(F+\na\X)dy,$$ where
$D$ is the unit disc in $\R^2.$ The formula above yields, by (\ref{6.5}),
\be W(QFR)=W(F)\qfq F\in M_+^{2\times2},\quad Q,\,\,R\in\SOII.\label{6.6}\ee

For $\var\in C^1(S,\R^3),$ consider a $C^1$ surface given by
$$S_\var=\{\,\var(x)\,|\,x\in S\,\}.$$ Denote the induced metric of $S_\var$ from $\R^3$ by
$g_\var=\<\cdot,\cdot\>\circ\var.$

We have the following.

\begin{pro}\label{p6.1}
Let the functional $I_0$ be given by Theorem $\ref{nt5.1}.$ Then
\be I_0(\var)=2\int_SW(d\var(x))dg\qfq\var\in C^1(S,\R^3),\label{6.10}\ee where
$$d\var(x)=\Big(\<E_i,\var_*e_j\>(\var(x))\Big)_{2\times2},$$ and, $e_1,$ $e_2$ and $E_1,$ $E_2$ are  positively orientated orthonormal bases of $S_x$ and  $(S_\var)_{\var(x)},$ respectively.
\end{pro}

{\bf Proof}\,\,\,
Let $E_3$ be the normal field of $S_\var.$ Then $\<E_3,\var_*e_i\>=0$ for $x\in S$ and $i=1,$ $2.$
Since $\var_*e_i=\sum_{j=1}^3\<E_j,\var_*e_i\>(\var(x))E_j$ for $1\leq j\leq3,$ i.e.,
$$ \Big(\,\var_*e_1\,|\,\var_*e_2\,|\,z\,\Big)=\Big(\,E_1\,|\,E_2\,|\,E_3\Big)\left(\begin{array}{ccc}\<E_1,\var_*e_1\>&\<E_1,\var_*e_2\>&\<E_1,z\>\\
\<E_2,\var_*e_1\>&\<E_2,\var_*e_2\>&\<E_2,z\>\\
0&0&\<E_3,z\>\end{array}\right),$$ for all $z\in\R^3,$ we obtain, by (\ref{6.1}),
$$\hat{W}_0(\var_*e_1,\var_*e_2)=W_0(d\var)\qfq x\in S.$$
The proof is complete. \hfill$\Box$\\

Next, let us assume that $M\subset \R^3$ is a $C^2$ surface with the induced metric $g.$ Suppose a middle surface of a shell $S$ is a bounded, open set of $M.$ We assume that all deformations of
$S$ are confined in $M.$ For such a deformation $\u,$ we assume that the stored energy is given
\be E(\u)=\int_SW(d\u)dg,\label{6.11}\ee where $W:$ $M_+^{2\times2}\rw\R$ satisfies (\ref{6.6}).

{\bf Membrane Shells of Revolution}\,\,\,Let $\psi$ be a $C^2$ function on $[0,\infty)$ with $\psi'(0)=0.$ Consider a surface of revolution given by
$$M=\{\,(x,\psi(r))\in\R^3\,|\,x=(x_1,x_2)\in\R^2,\,\,r=|x|\,\}.$$
The Gaussian curvature is
$$\kappa(p)=\frac{\psi'(r)\psi''(r)}{r(1+\psi'^2(r))^2}\qfq p=(x,\psi(r))\in M.$$ The normal field is
$$N(p)=\frac{1}{\sqrt{1+\psi'^2(r)}}(-\frac{\psi'(r)}{r}x,1).$$

Let $o=(0,\psi(0))\in M$ be fixed. Then $M_o=\R^2.$
Let $\zeta(t)$ be defined by the equation
\be t=\int_0^{\zeta(t)}\sqrt{1+\psi'^2(s)}ds\qfq t\geq0.\label{6.11}\ee Let
$$\gamma(t)=\Big(\zeta(t)v,\psi(\zeta(t))\Big)\qfq t\in\R,$$ where $v=(v_1,v_2)\in\R^2$ with $v_1^2+v_2^2=1.$

\begin{lem}\label{l6.1}
 $\gamma(t)$ is a normal geodesic such that
$$\gamma(0)=o,\quad\dot{\gamma}(0)=v.$$
\end{lem}

{\bf Proof}\,\,\,Let $D$ denote the connection of the induced metric $g$ of surface $M.$ Then
$$D_{\dot{\gamma}(t)}\dot{\gamma}=\ddot{\gamma}(t)-\<\ddot{\gamma}(t),N(\gamma(t))\>N(\gamma(t))=0\qfq t\geq0,$$ which prove the lemma.  \hfill$\Box$\\

It follows from Lemma \ref{l6.1} that
\be\kappa(t)=\kappa(\gamma(t))=\frac{\psi'(\zeta(t))\psi''(\zeta(t))}{\zeta(t)(1+\psi'^2(\zeta(t)))^2}\qfq t\geq0,\label{6.12}\ee where $\zeta$ is given by (\ref{6.11}).

In addition, we have the following.

\begin{pro} $(M,g,o)$ is a model.
\end{pro}

{\bf Proof}\,\,\,Since $n=2,$  we have
$${\bf R}(\dot{\gamma}(t),X,\dot{\gamma}(t),X)=\kappa(\gamma(t))|X|^2\quad\mbox{for all}\quad X  \in M_{\gamma(t)},\quad \<X,\dot{\gamma}(t)\>=0,$$ where ${\bf R}(\cdot,\cdot,\cdot,\cdot)$ is
the curvature tensor, which imply that formula (\ref{nn2.11}) holds true. By Proposition \ref{np2.3}, the proof is complete.  \hfill$\Box$\\

Let a middle surface $\B$ of a membrane shell be the unit geodesic disc in $(M, g)$ centered at the point $o,$ i.e.,
$$\B=\Big\{\,\Big(\zeta(t)v,\psi(\zeta(t))\Big)\,\Big|\,v\in\R^2,\,\,0\leq t<1\,\Big\}.$$ The radial deformations are given by
$$\u(p)=\Big(\zeta\circ\var(\rho)v,\,\,\psi\circ\zeta\circ\var(\rho)\Big)\qfq p=\Big(\zeta(\rho)v,\psi(\zeta(\rho))\Big),$$ where $\var$ is a function on $[0,\infty).$
Thus, all the theorems, corollaries, and propositions in Sections 2-4 hold true for radial deformations $\u.$ We do not repeat them here.

To end this section, we present two examples which verify the assumptions on the radial curvature in Proposition \ref{inp3.2} and Theorem \ref{t4.4}, respectively.

\begin{exl} Let $\varepsilon>0$ be  given and let $\psi_0\in C_0^\infty(0,\infty)$ be such that
$$\psi_0(t)=0\qfq 0\leq t\leq\varepsilon;\quad\psi_0(t)=1\qfq t\geq2\varepsilon.$$
Let
$$\psi_a(t)=\kappa_0+\frac{a\psi_0(t)}{1+t}\qfq t\geq0,$$ where $\kappa_0<0$ is a constant and $a>0$ is such that
$$\int_0^\infty s(\kappa_a)_+(s)ds=\int_\epsilon^{2\varepsilon}s(\kappa_a)_+(s)ds\leq1.$$ Consider the incompressible case. Let $\Phi$ satisfy the Baker-Ericksen inequalities with $\hat{\Phi}''(1)>0$ and such that $(\ref{cn3.5})$ is true. Let
$$I(A)=\int_\B W(d\u)dg-P\int_{\S}\rho(\u)d\S.$$ It follows from Proposition $\ref{inp3.2}$ $(i)$ that, for $A>0$ small and $P=\chi(A),$ $A$ is a local minimizer of $I(A).$
\end{exl}

\begin{exl} Let
$$\psi(s)=a\log(1+s^2)\qfq s\geq0,$$ where $a>0$ is a constant. If $0<a\leq1/\sqrt{2},$ then
\be\int_0^\infty t\kappa_+(t)dt\leq1,\quad\int_0^\infty t\kappa_-(t)dt<\infty.\label{6.14}\ee

By $(\ref{6.12}),$ we have
$$\kappa(t)=4a^2\frac{1-\zeta^4(t)}{[1+(4a^2+2)\zeta^2(t)+\zeta^4(t)]^2}\qfq t\geq0.$$
Let $t_0>0$ be given by $\zeta(t_0)=1.$ Thus,
\beq \int_0^\infty t\kappa_+(t)dt&&=4a^2\int_0^{t_0}\frac{t[1-\zeta^4(t)]}{[1+(4a^2+2)\zeta^2(t)+\zeta^4(t)]^2}dt\nonumber\\
&&\leq 4a^2\int_0^{1}\frac{1-\zeta^2}{1+(4a^2+2)\zeta^2+\zeta^4}d\zeta\leq 4a^2\int_0^1\frac{1-\zeta^2}{(1+\zeta^2)^2}d\zeta=2a^2.\nonumber\eeq
Similar arguments yield the second estimate in $(\ref{6.14}).$

Consider the compressible case. Let $W$ be given by $(\ref{3.2})$ with $n=2$ such that $(A 1)$-$(A 9)$ hold true.
It follows from Theorem $\ref{t4.4}$ that, for any $\lam$ sufficiently large, a minimizer of $I$ on $\A_\lam$ is a cavitating equilibrium solution, where $I$ is given by $(\ref{4.1}).$
\end{exl}

\setcounter{equation}{0}
\section{Cavitation for Ellipsoids in $\R^n$}
\def\theequation{6.\arabic{equation}}
\hskip\parindent Let $$g(x)=G(x)$$ be a symmetric, $C^2,$ and positively definite matrix for each $x\in\R^n$
and regard the pair $(\R^n,g)$ as a Riemannian manifold. The  study in the previous sections describes radial deformations of a ball-like body  if we apply it to the Rimannian manifold $(\R^n,g).$ The existence of corresponding cavitating equilibrium solutions depends on the geometric properties
of the metric $g$ and on the growth properties of  the constitutive function $W$ together.

We denote the metric $g$ on $\R^n$ by
$$g(X,Y)=\<X,Y\>_g=\<G(x)X,Y\>\qfq X,\,\,Y\in \R^n_x=\R^n,$$ where $\<\cdot,\cdot\>$ is the Euclidean metric in $\R^n.$ Then under the natural coordinates $x=(x_1,\cdots,x_n)$
$$g_{ij}(x)=\<\frac{\pl}{\pl x_i},\frac{\pl}{\pl x_j}\>_g=\<G(x)\frac{\pl}{\pl x_i},\frac{\pl}{\pl x_j}\>\qfq x\in\R^n.$$
Consider a body which occupies the open subset $\Om$ of $\R^n.$
 A map $\u:$ $\Om\rw \R^n$ is said to be a deformation of the body $\Om.$

\begin{thm}\label{t5.1} Let $W:$ $M_+^{n\times n}\rw\R$ be a constitutive function and satisfy $(\ref{2}).$ Then
\be W(d\u)=W(G^{1/2}(\u(x))\nabla\u(x)G^{-1/2}(x)),\label{5.1}\ee
where $\u=(u_1,\cdots,u_n)$ is a deformation, $W(d\u)$ is defined by $(\ref{3})$ on the Riemannian manifold $(\R^n,g)$ and $$\nabla u=\Big(\frac{\pl u_i}{\pl x_j}(x)\Big)\qfq x\in\R^n$$
is the gradient matrix of the map $\u$ in the Euclidean space $\R^n.$
\end{thm}

{\bf Proof}\,\,\,Let $x\in\R^n$ be given. Let $\{e_i\}$ and $\{E_i\}$ be orthonormal
bases of $(\R^n_x,g(x))$ and $(\R^n_{\u(x)},g\circ\u(x))$ with
positive orientation, respectively.

Let $$e_i=\sum_{j=1}^n\a_{ij}\pl_{x_j}|_x,\quad E_i=\sum_{j=1}^n\b_{ij}\pl_{x_j}|_{\u(x)}.$$
Then the relations, $$\delta_{ij}=\sum_{k=1}^n\a_{ik}\<\pl_{x_k},e_j\>_g,\qfq 1\leq i,\,\,k\leq n,$$ imply that \be
\Big(\a_{ij}\Big)\Big(\<\pl_{x_i},e_j\>_g\Big)=I.\label{1.7}\ee
On the other hand, the relations,
$$\<e_i,\pl_{x_j}\>_g=\sum_{k=1}^n\a_{ik}\<\pl_{x_k},\pl_{x_j}\>_g=\sum_{k=1}^n\a_{ik}g_{kj}(x),\qfq1\leq i,\,\,j\leq n,$$
yield \be \Big(\<e_i,\pl_{x_j}\>_g\Big)=\Big(\a_{ij}\Big)G(x),\label{1.8}\ee
where $G(x)=\Big(g_{ij}(x)\Big)$.

It follows from  formulas (\ref{1.7}) and (\ref{1.8}) that $$\Big(\a_{ij}\Big)G(x)\Big(\a_{ij}\Big)^T=I,$$
where the superscript "T" denotes the transpose.
A similar computation gives
$$\Big(\b_{ij}\Big)G(\u(x))\Big(\b_{ij}\Big)^T=I.$$

Noting the relations $$\u_*\pl_{x_i}=\sum_{k=1}^nu_{kx_i}\pl_{x_k}|_{\u(x)},$$
we obtain $$d\u(E_i,e_j)=\<E_i,\u_*e_j\>_g\circ\u(x)=\sum_{klp}\b_{ik}g_{kp}(\u(x))u_{px_l}\a_{jl},$$ that is,
$$\Big(d\u(E_i,e_j)\Big)=\Big(\b_{ij}\Big)G(\u(x))\nabla\u(x)\Big(\a_{ij}\Big)^T=QG^{1/2}(\u(x))\nabla u(x)G^{-1/2}(x) R,$$
where $$Q=\Big(\b_{ij}\Big)G^{1/2}(\u(x)),\quad R=G^{1/2}(x)\Big(\a_{ij}\Big)^T,$$ belong to $\SO(n).$
Thus,   formula (\ref{5.1}) follows from  assumption
(\ref{2}).  \hfill$\Box$\\

{\bf  Total Stored Energy }\,\,\, Let $dg$ denote the volume element of $\R^n$ in the metric $g.$ Then
$$dg={\det}^{1/2}G(x)dx\qfq x\in\R^n$$ where $dx$ is the volume element of $\R^n$ in the Euclidean metric.
Let $W:$ $M_+^{n\times n}\rw\R$ be a constitutive function and satisfy (\ref{2}).
By Theorem \ref{t5.1}, in a typical deformation in which the particle $x\in \Om$ is
displaced to $\u(x)\in \R^n$ energy (\ref{4}) becomes \be
E(\u)=\int_\Om \hat{W}\Big(x,\u(x),\nabla\u(x)\Big)dx,\label{5.4}\ee
where
\be\hat{W}(x,y,F)=W(G^{1/2}(y)FG^{-1/2}(x)){\det}^{1/2}G(x),\label{5.5}\ee for $x,\,\,y\in\R^n$ and  $F\in M_+^{n\times n},$ is the total constitutive law.
Then
$$\hat{W}(x,y,F)=\Phi\Big(v_1(x,y,F),\cdots,v_n(x,y,F)\Big)\det^{1/2}G(x),$$ where $v_1(x,y,F),$ $\cdots,$ $v_n(x,y,F)$ denote the singular vaules of
$G^{1/2}(x)FG^{-1/2}(y)$ for $x,$ $y\in\R^n$ and $F\in M_+^{n\times n}.$
Thus, formula (\ref{5.4})  is composed by the constitutive function $W$
and the matrices  $G(x)$ together.

\begin{exl} Introduce a metric on $\R^n$ by
\be g=e^{2a(r)}I\qfq r=|x|\in\R^n,\ee where $a(s)$ is a $C^2$ function on $[0,\infty)$ and $I$ is the unit matrix.  Suppose a constitutive function
 $W=\Phi$ is given by $(\ref{3.2}).$
Then the energy ensity $(\ref{5.5})$ is
$$\hat{W}(x,y,F)=W\Big(q(x,y)F\Big)e^{a(|x|)}=\sum_{i=1}^n\phi\Big(q(x,y)v_i\Big)e^{na(|x|)}+h\Big(q^n(x,y)v_1\cdots v_n\Big)e^{na(|x|)}$$
for $x,$ $y\in\R^n,$
where $$q(x,y)=\frac{e^{a(|y|)}}{e^{a(|x|)}}$$ and $v_1,$ $\cdots,$ $v_n$ denote the singular vaules of $F$ for $F\in M_+^{n\times n}.$
\end{exl}

We take $o$ to be the origin $0$ in $\R^n$ to consider what conditions on $g$ are needed for $(\R^n,g,o)$ being a model.

\begin{pro} Let $G(x)$ satisfy
\be G(x)x=b(\eta(x))Ax\qfq x\in\R^n,\label{5.7}\ee where $\eta(x)=\sqrt{\<Ax,x\>},$ $A$ is a sysmetric, positive, and constant matrix, and $b$ is a positive $C^2$ function on $[0,\infty)$ with $b(0)=1.$ Then

$(i)$\,\,\, $(\R^n,g,o)$ is a model;

$(ii)$\,\,\,A geodesic ball of $(\R^n,g)$ centered at $o$ is an ellipsoids of the Euclidean space $\R^n.$
\end{pro}

{\bf Proof}\,\,\,(i)\,\,\, We introduce one more metric $g_1$ on $\R^n$ by
$$g_A(X,Y)=\<AX,Y\>\qfq X,\,\,Y\in\R^n_x,\,\,x\in\R^n.$$ Clearly, (\ref{5.7}) implies $G(0)=g_A.$
We now have three metrics on $\R^n,$ that are $g,$ $g_A,$ and the Eucliean metric $\<\cdot,\cdot\>.$

Since $A$ is constant, for $x\in\R^n$ given, $x\not=0,$  the curve
$$\a(t)=t\frac{x}{\eta(x)}\qfq t\geq0,$$ is a normal geodesic in $(\R^n,g_A)$ initiating at $o.$
Thus, $\eta(x)$ is the distance function of $(\R^n,g_A)$ from $x\in\R^n$ to $o$ and
$$\na_A\eta=\frac{x}{\eta}\qfq x\in\R^n,\,\,x\not=0,$$ where $\na_A$ denotes the Levi-Civita connection of $(\R^n,g_A).$

Let $\si$ be the solution to problem
$$\si'(t)=b^{-1/2}(\si(t))\qfq t>0,\quad \si(0)=0.$$
Set
\be\gamma(t)=\si(t)X(x)\qfq t>0,\label{5.8}\ee where $X(x)=\na_A\eta.$ Next, we prove that $\gamma(t)$ is a normal geodesic of $(\R^n,g).$

Let $D$ be the Levi-Civita connection of $(\R^n,g).$ We compute $D_XX.$ Let $x_0\in\R^n$ be given.  Let $Z$ be a constant vector satisfying $\<x_0,Z\>_{g(x_0)}=0.$ Then
$$\<(\na_A\eta)(x_0),Z\>_A=\frac{1}{\eta}\<Ax_0,Z\>=\frac{1}{\eta(x_0)b(\eta(x_0))}\<G(x_0)x_0,Z\>=0,$$
$$[X(x_0),Z]=(\na_A)_{X(x_0)}Z-(\na_A)_Z\na_A\eta=-\frac{1}{\eta(x_0)}Z\quad \mbox{( since
$A$ is constant)},$$
$$\<X,Z\>_g=\<G(x)\frac{x}{\eta},Z\>=b(\eta)\<\na_A\eta,Z\>_A.$$
We have
\beq\<D_XX,Z\>_{g(x_0)}&&=X\<X,Z\>_g-\<X,D_XZ\>_g=X(b)\<\na_A\eta(x_0),Z\>_A-\<X,D_ZX\>_g-\<X,[X,Z]\>_g\nonumber\\
&&=-\frac{1}{2}Z|X|_g^2=-\frac{1}{2}Z(b|\na_A\eta|_A^2)=-\frac{b'(\eta)}{2}\<\na_A\eta,Z\>_A=0,\nonumber\eeq for all $Z\in\R^n$ satisfying
$\<x_0,Z\>_{g(x_0)}=0.$ Thus, we obtain
$$D_{X(x_0)}X=\<D_XX,\frac{x_0}{|x_0|_{g(x_0)}}\>_{g(x_0)}\frac{x_0}{|x_0|_{g(x_0)}}=\frac{b'(\eta)}{2b(\eta)}X(x_0).$$
It follows that
$$D_{\dot{\gamma}_0(t)}\dot{\gamma}_0=\ddot{\si}(t)X+\dot{\si}^2(t)D_{X(x_0)}X=[\ddot{\si}(t)+\frac{b'(\eta)}{2b(\eta)}]X(x_0)=0,$$ where
$\gamma_0(t)=\si(t)\frac{x_0}{\eta(x_0)}.$

Let $\rho$ be the distance function from $x\in\R^n$ to $o$ in the metric $g.$ By (\ref{5.8}), we obtain
$$\rho(x)=\si^{-1}(\eta(x))\qfq x\in\R^n,$$ and
\be\exp_ot\frac{x}{\eta(x)}=\si(t)\frac{x}{\eta(x)}\qfq t\geq0,\,\,x\in\R^n.\label{5.9}\ee

Let $\psi:$ $(\R^n,g_A)\rw(\R^n,g_A)$ be a linear isometry. Then the operator $\Psi$ of (\ref{d2.11}) is given by
$$\Psi(x)=\psi x\qfq x\in\R^n,$$ which imply that $\Psi:$
$\exp_o\Sigma(o)\rw\exp_o\Sigma(o)$ is an isometry.

(ii)\,\,\,Let $\B(t)$ be the geodesic ball centered at $o$ with radius $t>0.$ It follows from (\ref{5.9}) that
\beq\B(t)&&=\Big\{\,\si(s)\frac{x}{\eta(x)}\,\Big|\,x\in\R^n,\,0\leq s<t\,\Big\}\nonumber\\
&&=\Big\{\,y\,\Big|\,y\in\R^n,\,\sqrt{\<Ay,y\>}<\si(t),\,0\leq s<t\,\Big\},\label{5.10}\eeq
that are ellipsoids  in $\R^n.$  \hfill$\Box$

\begin{rem} If $A=\diag\{\frac{1}{a_1^2},\cdots,\frac{1}{a_n^2}\},$ then the geodesic balls $(\ref{5.10})$ are
$$\B(t)=\Big\{\,y\,\Big|\,y\in\R^n,\,\sum_{i=1}^n\frac{y_i^2}{a_i^2}<\si^2(t)\,\Big\}$$ and the geodesic spheres are
$$\Big\{\,y\,\Big|\,y\in\R^n,\,\sum_{i=1}^n\frac{y_i^2}{a_i^2}=\si^2(t)\,\Big\}\qfq t>0.$$
\end{rem}

\begin{rem}Let
$$\Ga_{ij}(x)=x_i\pl x_j-x_j\pl x_i\quad x\in\R^n,\,\,1\leq i,\,j\leq n.$$ Let $A$ be a constant matrix. Then matrices
$$G(x)=b_1(\eta)A+b_2(\eta)Ax\otimes Ax+\sum_{ij}b_{ij}(x)\Ga_{ij}(x)\otimes\Ga_{ij}(x)\qfq x\in \R^n$$ meet conditions
$(\ref{5.7})$ since
$$G(x)x=[b_1(\eta)+b_2(\eta)\eta^2]Ax\qfq x\in\R^n.$$
\end{rem}

Let $G(x)$ satisfy assumptions (\ref{5.7}). Then radial deformations are given by
$$\u(x)=\si\circ\var(\rho)\frac{x}{\eta(x)}\quad \rho=\si^{-1}(\eta(x)),\,\,x\in\R^n.$$
Thus, all the theorems, corollaries, and propositions  in Sections 2-4 hold true for radial deformations above if we apply them to the model $(\R^n, g, o).$

Finally, let us consider some situations for which the radial curvature assumptions in Theorem \ref{t4.4} hold. Let $\kappa$ be a $C^1$ function
in $[0,\infty)$ such that
$$ \int_0^\infty s\kappa_+(s)ds\leq 1,\quad \int_0^\infty s\kappa_-(s)ds<\infty,$$ where $\kappa_+=\max\{0,\kappa\}$ and $\kappa_-=\min\{0,-\kappa\}.$ Suppose $f$ is the solution to problem
$$f''(t)+\kappa(t)f(t)=0\qfq t>0;\quad f(0)=0,\quad f'(0)=1.$$

By similar arguments as in the proof of Proposition 4.2 in \cite{GrWu}, we obtain the following.

\begin{pro} Let $A$ be a symmetric, positive, and constant matrix. Let
$$G(x)=\frac{1}{f^2(\eta)}A+\frac{1}{\eta^2}\Big[1-\frac{f^2(\eta)}{\eta^2}\Big]Ax\otimes Ax\qfq x\in\R^n,$$ where $\eta=\sqrt{\<Ax,x\>}.$ Then

$(i)$\,\,\,$G(x)$ are symmetric and positive for all $x\in\R^n.$

$(ii)$\,\,\,$(\R^n,g,o)$ is a model.

$(iii)$\,\,\,The radial curvature is $\kappa(t)$ for $t\geq0.$
\end{pro}

{}

\end{document}